\documentclass[11pt]{article}
\usepackage{amsgen,amsmath,amstext,amsbsy,amsopn, amssymb,amsfonts}
\usepackage{pstricks}
\usepackage{epsf}
\usepackage{float}
\usepackage{epsfig}
\usepackage{here}
\usepackage{fancyheadings}

\usepackage{fullpage}
\usepackage{amsfonts}
\usepackage{amssymb}

\newcommand{\iref}[1]{\eqref{#1}}

\def \be{\begin{equation}}
\def \ee{\end{equation}}
\newcommand{\beq}{\begin{eqnarray}}
\newcommand{\eeq}{\end{eqnarray}}
\newcommand{\bea}{\begin{array}{c}}
\newcommand{\eea}{\end{array}}
\newcommand{\bi}{\begin{itemize}}
\newcommand{\ei}{\end{itemize}}

\newtheorem{rem}{Remark}

\newtheorem{lem}{Lemma}

\newtheorem{theo}{Theorem}

\def\N{\mathbb N}
\def\R{\mathbb R}

\def\C{\mathbb C}

\def\vp{\varphi}
\def\t{\tilde}
\def\e{\varepsilon}
\def\eps{\varepsilon}
\def\ep{\epsilon}

\def\Chi{\raise .3ex \hbox{\large $\chi$}} 
\def\vp{\varphi}

\def\b{\beta}

\def\p{\partial}

\def\ii{\infty}

\def\BA{\mathcal{A}}

\def\BR{\mathcal{R}}
\def\BL{\mathcal{L}}

\def\BD{\mathcal{D}}

\def\BS{\mathcal{S}}
\def\BC{\mathcal{C}}
\def\BO{\mathcal{O}}

\def\b{\ \ \ \ $\Box$}

\title{Geometric optics  and boundary layers \\ for  Nonlinear-Schr\"odinger Equations.}
\author{D. Chiron, \,\,\, F. Rousset\footnote{Laboratoire J.A. DIEUDONNE, Universit{\'e} de Nice - Sophia Antipolis,
Parc Valrose, 06108 Nice Cedex 02, France, chiron@unice.fr, frousset@unice.fr}}
\date{}

\begin{document}

\maketitle
\begin{abstract}
We justify supercritical geometric optics in small time for the defocusing 
semiclassical Nonlinear Schr\"odinger Equation for a large class of 
non-necessarily homogeneous nonlinearities. The case of a half-space with 
Neumann boundary condition is also studied.
\end{abstract}

\section{Introduction}
\ \indent We consider the nonlinear Schr\"odinger equation in $\Omega \subset \R^d$
\be
\label{NLS}
i \e \frac{\p \Psi^\e}{\p t} + \frac{\e^2}{2} \Delta \Psi^\e 
- \Psi^\e f(|\Psi^\e|^2) = 0,
\indent \indent \Psi^\e : \R^+ \times \Omega \to \C
\ee
with an highly oscillating  initial datum under  the form
\be
\label{NLSCI}
\Psi^\e_{|t=0} = \Psi^\e_0 = a_{0}^\eps  \exp\bigg( \frac{i}{\e} \vp_0^\eps \bigg),
\ee
where  $\vp_0^\eps$  is  real-valued. We are interested 
in the semiclassical limit $\e \to 0$. The nonlinear Schr\"odinger equation \iref{NLS} 
appears, for instance, in optics, and also as a model for Bose-Einstein condensates, 
with $f(\rho) = \rho-1$, and the equation is termed Gross-Pitaevskii equation, 
or  also  with $f(\rho)= \rho^2$ (see \cite{KNSQ}). Some more complicated
 nonlinearities are also used especially in low dimensions, see \cite{KL}.

At first, let us focus on the case $\Omega= \mathbb{R}^d$. To guess the formal limit, 
when $\eps$ goes to zero, it is classical to use the {\it Madelung transform}, i.e to 
seek for a solution of \eqref{NLS} under the form
$$ \Psi^\e = \sqrt{\rho^\e} \exp\bigg( \frac{i}{\e} \vp^\e\bigg).$$
By separating real and imaginary parts an by introducing $u^\eps \equiv \nabla \varphi^\eps$, 
this allows to rewrite \iref{NLS} as an hydrodynamical system
\be
\label{Madelungrho}
\left\{\begin{array}{ll}
\displaystyle{\p_t \rho^\e + \nabla \cdot \big(\rho^\e u^\e\big)}  = 0 \\ \ \\
\displaystyle{ \p_t u^\e+ \big( u^\e \cdot \nabla \big) u^\e + \nabla \big( f (\rho^\e) \big)}  = \displaystyle{\frac{\e^2}{2}} \nabla \bigg( \displaystyle{\frac{\Delta \sqrt{\rho^\e}}{\sqrt{\rho^\e}}} \bigg).
\end{array}
\right.
\ee
The system \iref{Madelungrho} is a compressible Euler equation with an additional term in 
the right-hand side called {\it quantum pressure}. As $\e$ tends to $0$, the quantum pressure 
is formally  negligible  and \iref{Madelungrho} reduces to the (compressible) Euler equation
\be
\left\{\begin{array}{ll}
\label{Euler}
\p_t \rho + \nabla \cdot \big( \rho u \big)  = 0\\ \\
\p_t u+ \big( u \cdot \nabla \big) u + \nabla \big( f (\rho) \big)  =0. \\
\end{array}\right.
\ee
The justification of this formal computation has received much interest recently. 
The case of analytic data was solved in \cite{Gerard}. Then for data with Sobolev 
regularity and a defocusing nonlinearity, so that \iref{Euler} is hyperbolic, it was noticed by 
Grenier, \cite{G}, that it is more convenient to use the transformation
\be
\label{grenier}
\Psi^\eps= a^\eps \exp\big(i \frac{\varphi^\eps}{\eps} \big)
\ee
and to allow the amplitude $a^\eps$ to be complex. By using an identification between 
$\mathbb{C}$ and $\mathbb{R}^2$, this allows to rewrite \eqref{NLS} as
\be
\left\{\begin{array}{ll}
\label{eulergren}
\p_t a^\eps + u^\eps \cdot \nabla a^\eps + \displaystyle{\frac{a^\eps}{2}}\, \nabla \cdot u^\eps
 = \eps J \, \Delta a^\eps \\ \\
\p_t u^\eps+ \big( u^\eps \cdot \nabla \big) u^\eps + \nabla \big( f ( |a^\eps|^2) \big)  =0,
\end{array} \right.
\ee
where $J$ is the matrix of complex multiplication by $i$:
$$ J = \left(\begin{array}{cc} 0 & - 1 \\ 1 & 0 \end{array} \right).$$
When $\eps=0$, we find the system 
\be
\left\{\begin{array}{ll}
\label{Eulera}
\p_t a + u \cdot \nabla a + \displaystyle{\frac{a}{2}}\, \nabla \cdot u  = 0\\ \\
\p_t u+ \big( u \cdot \nabla \big) u + \nabla \big( f ( |a|^2) \big)  =0,
\end{array}\right.
\ee
which is another form of \eqref{Euler}, since then $(\rho \equiv |a|^2, u)$ solves \eqref{Euler}. 
The rigorous convergence of \eqref{eulergren} towards \eqref{Eulera} provided the initial 
conditions suitably converge was rigorously performed by Grenier \cite{G} in the case $f(\rho)= \rho$ 
(which corresponds to the cubic defocusing NLS). More precisely, it was proven in \cite{G} that 
there exists $T>0$ independent of $\eps$ such that the solution of \eqref{eulergren} is 
uniformly bounded in $H^s$ on $ [0,T]$. In terms of the unknown $\Psi^\eps$ of \eqref{NLS}, 
this gives that
$$ \sup_{ \eps \in (0, 1]}\sup_{[0, T]} \big|\!\big| \Psi^\eps \exp\big(- i \frac{\varphi}{\eps} \big) 
\big|\!\big|_{H^s } < + \infty$$
for every $s$ where $(a, u= \nabla \varphi$) is the solution of \eqref{Eulera}.
Furthermore, the justification of WKB expansions 
 under the form
 $$   \Psi^\eps -  \Big(\sum_{k=0}^m \eps^k a^k \Big)e^{i \varphi \over \eps}
  = \mathcal{O}(\eps^m)$$
   for every  $m$
was performed in \cite{G}. The main idea in the  work of Grenier \cite{G} is to use the symmetrizer 
$$ S \equiv \mbox{ diag }\Big( 1, 1, \frac{1}{4}  f'(|a|^2), \cdots, \frac{1}{4} f'(|a|^2) \Bigr)$$
of the hyperbolic system \eqref{Eulera} to get $H^s$ energy estimates which are uniform in $\eps$ 
for the singularly perturbed system \eqref{eulergren}. The case of nonlinearities for which 
$f'$ vanishes at zero (for instance the case $f(\rho)= \rho^2$) was left opened in \cite{G}. 
The additional difficulty is that for such nonlinearities, the system \eqref{Eulera} is only 
weakly hyperbolic at $a=0$ and in particular the symmetrizer $S$ is not anymore positive 
definite at $a=0$.


In more recent works, see \cite{Z}, \cite{LZ}, \cite{AC0} it was proven that for every 
weak solution of \eqref{NLS} with $f(\rho)=\rho-1$ or $f(\rho)= \rho$, the limits as 
$\e \rightarrow 0$
\beq
\label{weak} |\Psi^\eps |^2 - \rho \to 0 \quad \quad {\rm in} \ \ L^\ii([0,T],L^2) \quad \quad \quad \quad
\eps  {\it Im}\, \big( \bar{\Psi}^\e \nabla \Psi^\e \big) - \rho u \to 0 
\quad \quad {\rm in} \ \ L^\ii([0,T],L^1_{loc} )
\eeq
hold under some suitable assumption on the initial data. The approach used in these papers is 
completely different, and relies on the modulated energy method introduced in \cite{Brenier}. 
The advantage of this powerfull approach is that it allows to describe the limit of weak solutions 
and to handle general nonlinearities once the existence of a global weak solution in the energy 
space for \eqref{NLS} is known. Nevertheless, it does not give precise qualitative information 
on the solution of \eqref{NLS}, for example, it does not allow to prove that the solution remains 
smooth on an interval of time independent of $\eps$ if the initial data are smooth or to justify 
WKB expansion up to arbitrary orders in smooth norms.


In the  work \cite{AC}, the possibility of getting the same result as in \cite{G} for 
pure power nonlinearities $f(\rho)= \rho^\sigma$ in the case $\Omega= \mathbb{R}^d$ was studied. 
It was first noticed that, thanks to the result of \cite{MUK}, the system
\be
\label{aphi}
\left\{\begin{array}{ll}
\displaystyle{ \partial_{t} a + \nabla \varphi \cdot \nabla a + \frac{a}{2} \Delta\varphi} & = 0 \\
\displaystyle{ \partial_{t} \varphi +  \frac{1}{2}|\nabla \varphi |^2 + f( |a|^2) } & = 0,
\end{array}
\right.
\ee
with the initial condition $ \big( a , \vp \big)_{/t=0} = \big( a_{0},\varphi_{0}\big) \in H^\infty$ 
 has a unique smooth maximal solution $(a,\varphi) \in \BC\big( [0,T^*[,H^s(\R^d) \times H^{s-1}(\R^d) \big)$ for every $s$. It was then established:

\begin{theo}[\cite{AC}]
\label{theoAC}
Let $d \leq 3$, $\sigma \in \mathbb{N}^*$ and initial data $a_{0}^\eps$, $\varphi_{0}$ in 
$H^\infty$ such that, for some functions $(\vp_0,a_0) \in H^{\ii}$,
$$ \big|\!\big|a_{0}^\eps - a_{0} \big|\!\big|_{H^{s}} = \BO(\eps),$$
for every $s \geq 0$. Then, there exists $T^*>0$ such that \eqref{aphi} with $f(\rho) = \rho^\sigma$ 
has a smooth maximal solution $(a, \varphi) \in \BC([0,T^*[,H^\ii \times H^\ii)$. Moreover, 
there exists $T\in(0,T^*)$ independent of $\eps$, such that the solution of 
\eqref{NLS}, \eqref{NLSCI} remains smooth on $[0, T]$ and verifies the estimate
\be
\label{estAC}
\sup_{ \eps \in (0, 1] } \big|\!\big|\Psi^\eps 
\exp\big(-i \frac{\varphi}{\eps} \big) \big|\!\big|_{L^\infty([0, T], H^s )} < + \infty,
\ee
where
\begin{itemize}
\item if $\sigma = 1$, then $ s \in \mathbb{N}$ is arbitrary,
\item if $\sigma = 2$ and $d=1$, then one can take $s=2$,
\item if $\sigma = 2$ and $2 \leq d \leq 3$, then one can take $s=1$,
\item if $\sigma \geq 3$ then one can take $s= \sigma$.
\end{itemize}
\end{theo}

As emphasized in \cite{AC}, in some cases, the global existence of smooth solutions 
is already known for \eqref{NLS}. For example, in the quintic case, $\sigma =2$, global 
existence is known for $d \leq 3$ (see \cite{CKSTT} for the difficult critical case $d=3$), 
so that only the bound \eqref{estAC} is interesting. Nevertheless, Theorem \ref{theoAC} 
may be also applied to cases where \eqref{NLS} is $H^1$ super-critical 
($\sigma \geq 3$, $d=3$ for example) and hence the fact that it is possible to 
construct a smooth solution on a time interval independent of $\eps$ is already interesting. 
The main ingredient used in \cite{AC} is a subtle transformation of \eqref{NLS}
into a perturbation of a quasilinear  symmetric  hyperbolic system with non smooth 
coefficients when $\sigma \geq 2$.\\

The first aim of this paper is to prove that the estimate \eqref{estAC} holds true 
for every $s$, every dimension $d$ and every nonlinearity $f$ which satisfies the 
following assumption:
$$ (\BA) \indent \indent f\in \BC^\infty \big( [0,+\ii) \big), \quad  \quad \quad f(0)=0, 
\quad  \quad \quad f'>0 \ \ {\rm on} \ \ (0,+\ii), \quad \quad \quad 
\exists n \in \N^*, \ \ f^{(n)}(0) \not = 0.$$


Note that we allow $f'$ to vanish at the origin. The assumption $(\BA)$ takes into account in 
particular all the homogeneous polynomial nonlinearities $f(\rho) = \rho^\sigma$ but also 
nonlinearities under the form $f(\rho)= \rho^{\sigma_{1}} + \rho^{\sigma_{2}}$
 or $ {\rho^\sigma \over  1+ \rho}$ for example. Our result 
reads:

\begin{theo}
\label{Conver}
We assume $(\mathcal{A})$, and consider an initial data \eqref{NLSCI} with 
$a_{0}^\eps$, $\varphi_{0}^\eps$ in $H^\infty$ such that, for some real-valued functions 
$(\vp_0,a_0) \in H^{\ii}$, we have for every $s$,
$$ \big|\! \big| a^\e_0 - a_{0} \big|\! \big|_{H^{s}} = \BO(\e) \quad \quad {\it and} \quad \quad 
\big|\! \big| \varphi^\e_0 - \varphi_{0} \big|\! \big|_{H^{s}} = \mathcal{O}(\eps).$$

Then, there exists $T^*>0$ such that \eqref{Eulera} with initial value 
$(a_{0}, \varphi_{0})$ has a unique smooth maximal solution $(a, \varphi) \in \BC([0,T^*[,H^\ii \times H^\ii)$. 
Moreover, there exists $T\in (0,T^*]$ such that  for every $\eps \in (0, 1)$, the solution $\Psi^\e$ to \iref{NLS}-\iref{NLSCI} 
exists at least on $[0,T]$ and satisfies for every $s$
$$\sup_{ \eps \in (0, 1] } \bigg|\!\bigg|\Psi^\eps 
\exp\big(- \frac{i}{\eps} \varphi \big) \bigg|\!\bigg|_{L^\infty([0, T], H^s )} 
< + \infty.$$
More precisely, there exists $\varphi^\eps = \varphi + \mathcal{O}_{H^\infty}(\eps)$ such that, 
for every $s$,
\be
\label{estconv}
\bigg|\! \bigg| \Psi^\e \exp\big( -\frac{i}{\e} \vp^\e \big) - a 
\bigg|\! \bigg|_{L^{\ii}([0,T],H^s)} = \BO(\e).
\ee
\end{theo}

Let us give a few comments on the statement of Theorem \ref{Conver}.

At first, note that Theorem \ref{Conver} contains a result of local existence of 
smooth solutions for \eqref{aphi} in the case of non necessarily homogeneous nonlinearities 
satisfying $(\BA)$. Since $(a, \nabla \varphi$) solves a compressible type Euler equation, 
the case of a homogeneous nonlinearity was studied in \cite{MUK}, and we thus give an extension 
of this result to smooth non-linearities satisfying assumption $(\mathcal{A})$.
 A precise 
statement of our result with the required regularity of the initial data is given in 
Theorem \ref{theoeulernh} below. The new difficulty when $f$ is not homogeneous
 is that the nonlinear symmetrization does not  seem to allow to  transform
  the problem into a classical symmetric or symmetrizable hyperbolic system
   with smooth coefficients.

The correction of order $\eps$ that we have to add to the phase to get the estimate 
\eqref{estconv} is expected. Indeed, a perturbation of order $\eps$ in the phase modifies 
the amplitude at the leading order.


     Our approach to prove Theorem \ref{Conver} is completely different from the one of 
\cite{AC} and \cite{G}. We do not work any more on the system \eqref{eulergren} 
 or any reformulation of \eqref{NLS} into a perturbation of a quasilinear symmetric 
hyperbolic system, but directly on the NLS equation \eqref{NLS}. Basically, we first 
prove the linear stability for \eqref{NLS} in arbitrary Sobolev norms of highly oscillating 
solution of the form $a e^{i \varphi / \eps}$ and then use a fixed point argument to prove 
the nonlinear stability.  The crucial estimate of linear stability 
 of highly oscillating solution is given in Lemma \ref{dtN} and Theorem \ref{estiHS}.

 This actually allows to justify WKB expansions up to arbitrary orders 
(see Theorem \ref{WKBstab}).
Since we deal in this paper with sufficiently smooth and in particular bounded solutions, the 
assumption $(\BA)$ can be replaced by a local version where
 we  assume that $f'>0$  on $(0, \beta) $
 with   $\beta$ independent 
of $\eps$  if  the initial datum verifies $ | a_{0} |^2 < \beta$.
Indeed, since $a^0$ takes it values in the (weak) hyperbolic
 region   of the 
limit system \eqref{Eulera}, there still exists a  local smooth solution of \eqref{Eulera}
 defined on $[0, T]$ for some $T>0$ and the stability argument
  leading to Theorem \ref{Conver} still holds.
Consequently, our result can also be applied to nonlinearities like 
$f(\rho)= \rho^{\sigma_{1}}- \rho^{\sigma_{2}}$ for every $\sigma_{2}>\sigma_{1}$ provided 
$|a_{0}|^2 \leq \beta <1$. Note that when $\sigma_{2}$ is too large, the classical global 
existence result of weak solutions (see \cite{GV}) for \eqref{NLS} is not valid and hence it 
does not seem possible to use the modulated energy method of \cite{AC0}, \cite{LZ} to investigate the semi-classical limit.

\bigskip

Finally, the last advantage of our approach is that it can be easily generalized to 
the case of a domain with boundary and to non-zero condition at infinity. This will be the aim of the second part of the paper.
We shall restrict ourself to a physical case, the Gross-Pitaevskii equation, i.e. $f(\rho) = \rho -1$. 
The generalization to more general nonlinearities satisfying an assumption like $ (\mathcal{A})$ 
is rather straightforward. This simplifying assumption is only made to avoid the multiplication 
of difficulties. Again to avoid too many technicalities, we restrict ourselves to the simplest 
domain $\Omega= \mathbb{R}^d_{+}= \mathbb{R}^{d-1} \times (0, + \infty).$  For $x \in \R^d_+$, we 
shall use the notation $ x=(y, z), \, y \in \mathbb{R}^{d-1}, \, z >0$. We add to \eqref{NLS} 
the Neumann boundary condition
\beq
\label{dir+}
\partial_{z}\Psi^\eps(t, y, 0) = 0.
\eeq
We also impose the following condition at infinity
\beq
\label{inf+}
\Psi^\eps (t,x) \sim \exp 
\Big( - i \, t \, \frac{ |u^{\infty}|^2}{2\e} + i \, \frac{ u^\infty \cdot x}{\e} \Big), 
\quad \quad \quad |x| \rightarrow + \infty,
\eeq
that we can write in hydrodynamical variables
$$ \big| \Psi^\eps (t,x) \big|^2 \to 1, \quad \quad \quad u^\e(t,x) \to u^\ii, 
\quad \quad \quad |x| \rightarrow + \infty,$$
where $u^\infty$ is a constant vector. This condition appears naturally when we study a moving 
obstacle in the fluid. Indeed, if we start from \eqref{NLS} with the Neumann boundary condition 
on an obstacle moving at constant velocity and  fluid at rest at infinity, then we can use the 
Galilean invariance of \eqref{NLS} to transform the problem into the study of \eqref{NLS} in a 
fixed domain but with the condition \eqref{inf+} at infinity.

This problem with such boundary conditions is physically meaningfull
 since it can be used to describe superfluids past an obstacle
  (we refer to \cite{PNB} for example).
   The  semiclassical limit $\eps$ tends to zero was already studied
    in \cite{LZ} by using the modulated energy method.
     The limit \eqref{weak} was  proven  with $(\rho, u)$ the solution of the compressible
      Euler  equation with  boundary condition
       $ u\cdot n_{/\partial \Omega} = 0$, $n$ being the normal to the boundary.
         Note that the result of \cite{LZ} is restricted to the two-dimensional case
       only in order to  have a global solution in the energy space of
       \eqref{NLS}.  By using more recent results on the Cauchy problem,
        \cite{A},  one can also get the result  in the three-dimensional
         case  at least when $u^{\infty} = 0 $.
          Our aim here is to give a more precise description of the convergence
           which takes into account  boundary layers.
        More precisely, since the solution of the Euler system \eqref{aphi}
         cannot match the Neumann boundary condition
          $\partial_{z}a(t,y,0)=0$, a boundary layer of weak amplitude
           $\eps$ and of size $ \eps$ appears.  They are formally described
            for example in \cite{PNB}. WKB expansions $\Psi^\eps
            = a^\eps e^{ i \frac{\varphi^\eps}{\eps  } }$
            are thus to be seek under the form
\beq
\label{expintro}
a^\eps = a^0 +\sum_{k= 1}^m\eps^k \Bigl( a^k(t,x) + A^k(t,y, \frac{z}{\eps}) \Big), \quad 
\varphi^\eps= \varphi^0 + \sum_{k=1}^m \eps^k\Big(  \varphi^k(t,x) + \Phi^k(t,y, \frac{z}{\e}) \Big)
\eeq 
   where the profiles $A^k(t,y,Z)$, $\Phi^k(t,y,Z)$ 
    are exponentially decreasing in the $Z$ variable
     and are chosen such that
        $$   \partial_{z}a^k(t,y,0) + \partial_{Z}A^{k+1}(t,y,0) = 0, \quad
  \partial_{z}\varphi^k(t,y,0) + \partial_{Z}\Phi^{k+1}(t,y,0) = 0   $$ 
   so that the approximate WKB expansion $ \Psi^{WKB}=a^\eps \exp\big( \frac{i}{\e} \varphi^\e \big)$
    matches the Neumann boundary condition \eqref{dir+}.
     Our result  (Theorem \ref{theoGP+}) is that under suitable  assumptions on the initial conditions, we have
      the   nonlinear stability of WKB expansions:  in particular
       we have  the existence of a smooth solution for \eqref{NLS}, \eqref{dir+}, \eqref{inf+}
        on a time interval independent of $\eps$  and 
       the estimate
\beq
\label{grad+}
\big|\!\big| \Psi ^\eps e^{- i \frac{\vp^\e}{\e}} -a^\e \big|\!\big|_{W^{1, \infty}} \lesssim \eps.
\eeq
Note that it is necessary to incorporate the boundary layer $\e A^1$ in order to get 
\eqref{grad+} since its gradient has amplitude one in $L^\infty$. The case of Dirichlet boundary 
condition which is also  physically meaningfull, we again refer to \cite{PNB}, seems more 
complicated to handle as often in boundary layer theory in fluid mechanics since the boundary 
layers involved have amplitude one. This is left for future work.\\

The paper is organized as follows. In section \ref{lineaire}, we prove the linear 
stability in $H^s$ of an approximate WKB solution of \eqref{NLS} under the form 
$a^\eps \exp \big( i \frac{\varphi^\eps}{\eps} \big)$ in the case $\Omega = \mathbb{R}^d$. 
This is the crucial part towards the proof of Theorem \ref{Conver}. Next in section \ref{sectionWKB}, 
we give the construction of a WKB expansion up to arbitrary order and give the proof of the 
local existence of smooth solution for the compressible Euler equation with a pressure law 
satisfying $(\BA)$. In section \ref{nonlineaire}, we give the justification of WKB expansions at 
every order and recover Theorem \ref{Conver} as a particular case. This part uses in a 
classical way the linear stability result and a fixed point argument. Finally, in section 
\ref{half}, we study the problem in the half-space with Neumann boundary condition.


\section{Linear Stability}
\label{lineaire}

\ \indent In this section, we consider a smooth WKB approximate solution 
$\Psi^a = a^\eps \exp \big(i \frac{\varphi^\eps}{\eps} \big)$ of \eqref{NLS} such that
\beq
\label{resteapp}
NLS(\Psi^a)= R^\eps \exp \big(i \frac{\varphi^\eps}{\eps} \big),
\eeq
where
$$NLS(\Psi) \equiv i \eps \partial_{t} \Psi + \frac{\eps^2}{2} \Delta \Psi- \Psi f(|\Psi|^2).$$
Moreover, we also set
\begin{align}
\label{Rphi}
R_{\varphi} & \equiv \partial_{t} \varphi^\eps + \frac{1}{2}\, | \nabla \varphi^\eps |^2 +  f(|a^\eps |^2), \\ \label{Ra}
R_{a} & \equiv \partial_{t} a^\eps + \nabla \varphi^\eps \cdot \nabla a^\eps + \frac{1}{2}\, a^\eps \Delta \varphi^\eps,
\end{align}
so that
$$ R^\e = - a^\e R_\vp + i\e R_a.$$
Looking for an exact solution of \eqref{NLS} under the form
$$ \Psi^\eps = \Psi^a + w \, e^{ i \frac{\varphi^\eps}{\eps}} = (a^\eps + w) e^{i\frac{\varphi^\eps}{\eps}},$$ 
 we find that $w$ solves the nonlinear Schr\"odinger equation
\beq
\label{NLSw}
i \eps\Big( \p_t w + u^\eps \cdot \nabla w + \frac{1}{2}\, w\, \nabla \cdot u^\eps \Big) + 
\frac{\eps^2}{2} \Delta w - 2 (w,a^\eps)  f'(|a^\eps|^2) a^\eps = R_{\varphi} w - R^\eps + Q(w),
\eeq
where we have set
$$ u^\eps \equiv \nabla \varphi^\eps$$
and the nonlinear term $Q(w)$ is defined by
\beq
\label{Qw}
Q(w) \equiv (a^\eps + w ) \Big( f( |a^\eps + w|^2) - f(|a^\eps |^2) \Big) - 
2 (w,a^\eps) f'(|a^\eps|^2) a^\eps.
\eeq
Of course, $R^\eps$ will be very small and $R_{\varphi}$ (and $R_{a}$) are to be thought 
small (at least $\mathcal{O}(\eps)$) for applications to nonlinear stability results. 
Nevertheless, in this section the exact form of these terms is not important. The way 
to construct an accurate WKB solution will be explained in the next section.

\begin{rem}\rm
If we work with a non-linearity $f$ such that
 $f(|A|^2)=0$, we can impose a  non-zero condition at infinity  such as $a_0 \in A+ H^\ii$ and 
$\nabla \vp_0 \in U^{\infty} + H^\ii$. Since  we can still look for  the perturbation $w$ 
in $H^s$,  this does 
not affect the proofs.
\end{rem}

Since we expect the correction term $w$ to be small, we shall only consider in this 
section the linearized equation
\be
\label{NLSwlin}
i\e \frac{\p w}{\p t} + \BL^\e\,w = R_{\varphi} w +  F^\e, \quad x \in \mathbb{R}^d, 
\ee
where the linear operator $\mathcal{L}^\e$ is defined as
$$ \BL^\e(w) \equiv \frac{\e^2}{2} \Delta w + i\e\, u^\e \cdot \nabla w + 
\frac{i\e}{2}\, w \, \nabla \cdot u^\e - 2f'(|a^\e|^2) (w,a^\e) a^\e$$
and $(\cdot ,\cdot)$ stands for the real scalar product in $\C\simeq \R^2$.
In this section, $F^\eps$ is considered as a given source term. Of course, for the 
proof of Theorem \ref{Conver}, we shall apply the result of this section to
\be
\label{Feps}
F^\eps = - R^\eps + Q(w).
\ee
We notice that the first and last terms in $\BL^\e w$ are formally self-adjoint, and that the 
quadratic form (in $H^1$) associated to the operator
$$ \BS^\e\,w \equiv - \frac{\e^2}{2} \Delta w + 2 f'(|a^\e|^2) (w,a^\e) a^\e$$
is
$$\int_{\R^d} \big( w , \BS^\e\, w \big) = \frac{1}{2} 
\int_{\R^d} \e^2 |\nabla w|^2 + 4 f'(|a^\e|^2) (w,a^\e)^2.$$
It is then natural to consider the (squared) norm 
$\displaystyle{\int_{\R^d} \big( w , \BS^\e(w) \big)}$
 as a good energy for the linearized equation \iref{NLSwlin}. Consequently, we introduce
  the weighted norm
$$ N^\eps(w) \equiv \frac{1}{2} \int_{\R^d} \e^2 |\nabla w|^2 + 4 f'(|a^\e|^2) (w,a^\e)^2 +  K\e^2 |w|^2$$
for every  $K>0$ ($K$ will be chosen sufficiently large only in 
 the next subsection).

 Our first result of this section 
is a linear stability result in the energy norm $N^\eps(w)$.

\begin{lem}
\label{dtN}
 Assume that  $u^\e : [0,T] \times \R^d \to \R^d$ and $a^\e : [0,T] \times \R^d \to \C$ 
  are smooth and 
such that
$$ M \equiv \big|\!\big| \nabla_x u^\e \big|\!\big|_{L^\ii([0,T]\times \R^d)} 
+ \big|\!\big| \nabla_x ( \nabla \cdot u^\e ) \big|\!\big|_{L^\ii([0,T]\times \R^d)} 
+ \big|\!\big| |a^\e|^2 \big|\!\big|_{L^\ii([0,T]\times \R^d)} < +\ii.$$
Let $w\in \BC^1([0,T],H^2)$ be a solution of \iref{NLSwlin}. Then, there exists $C_M$ 
depending only on $d$, $f$ and $M$ such that
 for every $\eps \in (0, 1]$, the solution of \eqref{NLSwlin}
  satisfies the energy estimate
\begin{align}
\label{lin0}
\frac{d}{dt} N^\eps\big( w(t) \big) \leq & \ 
C_{M} \bigg( 1+ \frac{1}{\e} \big|\!\big| R_{a}(t) \big|\!\big|_{L^\ii} + 
\frac{1}{\e} \big|\!\big| R_{\varphi}(t) \big|\!\big|_{W^{1,\ii}} +
\frac{1}{\e^2} \big|\!\big| R_{\varphi}(t) \big|\!\big|_{L^\ii} \bigg) N^\eps\big( w(t) \big) \\
\nonumber & + \big|\!\big| F^\e(t) \big|\!\big|_{L^2}^2 - 
\int_{\R^d} \frac{4}{\e} f'(|a^\e|^2) (w,a^\e) (a^\e, i F^\e) + \int_{\R^d} ( \eps \Delta w, i F^\eps) .
\end{align}
\end{lem}

Note that it is very easy to get from \eqref{lin0} and the  Gronwall inequality a 
classical estimate of linear stability. Indeed, assuming  that
$R_a = \BO_{L^\ii([0,T],L^\ii)}(\e)$ and $R_\vp = \BO_{L^\ii([0,T],W^{1,\ii})}(\e^2)$
(which is true if $(a^\eps, \varphi^\eps)$  come from the WKB method), we infer from 
a crude estimate for the two last terms in \eqref{lin0} that for $0\leq t\leq T$,
$$\frac{d}{dt} N^\eps\big( w(t) \big) \leq 
C N^\eps\big( w(t) \big) + \frac{1}{\eps^2} \big|\!\big| F^\e(t) \big|\!\big|_{H^1}^2,$$
which gives for $0\leq t\leq T$
$$ N^\eps\big(w(t)\big) \leq e^{Ct}\Big( N^\eps\big(w(0)\big) + 
\frac{1}{\e^2} \int_{0}^t \big|\!\big|F^\e(\tau)\big|\!\big|_{H^1}^2 \ d\tau\Big),$$
which is a  more classical result of linear stability in the energy norm $N^\e(w)$ since 
the amplification rate $C$ is independent 
of $\eps$.
Nevertheless, to 
get $H^s$ estimates and  the best nonlinear results as possible, it is important to have 
the special structure of the two last terms in \eqref{lin0}.

Modulated linearized functionals like $N^\eps$ were also used in asymptotic
problems in fluid mechanics, see \cite{GrenierFM} for example.

\subsection{Proof of Lemma \ref{dtN}}

The norms $L^\ii$, $W^{1,\ii}$, $L^2$ ...   always stand for the  norms in  the $x$ variable. 

At first, 
since $\mathcal{S}^\eps$ is self adjoint, we have
\beq
\label{estdtN0}
\frac{d}{dt} \int_{\R^d} \big( \mathcal{S}^\eps w, w \bigr) = 
\int_{\R^d} 2 \big( \mathcal{S}^\eps w,\partial_{t}w \big) + 
2 \partial_{t} \big[f'(|a^\eps|^2)\big] (w, a^\eps)^2 + 
4f'(|a^\eps|^2) (w, a^\eps) (w, \partial_{t} a^\eps) .
\eeq
Next, we use \eqref{NLSwlin} to express $\partial_{t}w$ as
$$ \partial_{t} w = - \frac{i}{\eps} \mathcal{S}^\eps w  - 
\big( u ^\eps \cdot \nabla w + \frac{1}{2}\, w\, \nabla \cdot u^\eps \big) - 
\frac{ i}{\eps} R_{\varphi} w  - \frac{i}{\eps }F^\eps$$
to get
\beq
\label{estdtN1}
2\int_{\R^d} \big(\mathcal{S}^\eps w ,\partial_{t}w \big) 
 = 2 \int_{\R^d} \Big( \frac{\eps^2}{2} \Delta w - 2 f'(|a^\eps|^2) (w,a^\e) a^\e, u^\eps \cdot \nabla w 
+ \frac{1}{2}\, w \, \nabla \cdot u^\eps + \frac{i}{\eps } R_{\varphi} w + \frac{i}{\eps } F^\eps \Big).
\eeq
We shall now estimate the various terms in the right-hand side of \eqref{estdtN1}. Integrating 
by parts, we get
\begin{eqnarray*}
\int_{\mathbb{R}^d} \big( \eps^2 \, \Delta w, \frac{i}{\e} R_{\varphi} w \big) 
& = & - \e \int_{\mathbb{R}^d} \big( \nabla w, i w \nabla R_{\varphi} \big) \\
& \leq & \e\, |\!| \nabla R_{\varphi} |\!|_{L^\infty} |\!| w |\!|_{L^2}  |\!| \nabla w |\!|_{L^2} \\
& \leq & \frac{1}{\e}\, |\!| R_{\varphi}|\!|_{W^{1,\ii}} N^\eps(w).
  \end{eqnarray*}
Note that  we have used that $R_{\varphi}$ is real-valued
 and thus that
 $$ (\nabla w , i R_{\varphi} \nabla w ) = 0$$
  for the first equality. We also 
easily obtain  by integration by parts that 
\begin{eqnarray*}
\int_{\mathbb{R}^d} \Big(\eps^2 \Delta w , w\, \nabla \cdot u^\eps \Big) 
& \leq & C \Big( \big|\!\big| \nabla \cdot u^\eps \big|\!\big|_{L^\infty} + 
\big|\!\big| \nabla ( \nabla \cdot u^\eps) \big|\!\big|_{L^\infty} \Big)
 \, \Big( \eps^2 \big|\!\big|\nabla w \big|\!\big|_{L^2}^2 + \eps^2 \big|\!\big|w\big|\!\big|_{L^2}^2 \Big) \\
& \leq & C_M N^\eps (w).
\end{eqnarray*}
In the proof, $C_M$ is a harmless number which changes from line to line and which depends 
only on $M$. In particular, it is independent of $\eps$. Moreover, 
we can also write
for  $k= 1, \cdots, d$, 
\begin{eqnarray*} 
\int_{\mathbb{R}^d}  \bigl( \partial_{kk}^2  w , u^\eps \cdot \nabla w \bigr)
& =& -   \int_{\mathbb{R}^d} u^\eps \cdot \nabla  { |\partial_{k} w|^2 \over 2}
   - \int \bigl(\partial_{k}w, \partial_{k}u^\eps\cdot \nabla w \bigr) \\
   & =&   \int_{\mathbb{R}^d}  { |\partial_{k} w|^2 \over 2} \, \nabla \cdot u^\eps  
 - \int_{\mathbb{R}^d} \bigl(\partial_{k}w, \partial_{k}u^\eps\cdot \nabla w \bigr)
 \end{eqnarray*}
and hence,  we immediately infer
$$ \int_{\mathbb{R}^d} \Big( \eps^2 \Delta w , u^\eps \cdot \nabla w \Big) \leq C_M N^\eps (w).$$
Furthermore, from the inequality $2ab \leq a^2 + b^2$, there holds 
\begin{eqnarray}
- \frac{4}{\e} \int_{\mathbb{R}^d} f'(|a^\e|^2) (w,a^\e) \big( a^\e, i R_{\varphi} w \big) 
& \leq & \frac{C_M}{\e^2} |\!|R_{\varphi}|\!|_{L^\infty} 
\int_{\mathbb{R}^d} \big( f'(|a^\eps|^2) \big)^{\frac{1}{2} } \big| (w, a^\eps ) \big|\, \eps | w | \nonumber \\
 & \leq & \frac{C_M}{\e^2} |\!|R_{\varphi}|\!|_{L^\infty} 
\int_{\mathbb{R}^d } f'(|a^\eps|^2) (w, a^\eps)^2 + \e^2 | w|^2 \nonumber \\
\label{estdtN2} & \leq & \frac{C_M}{\e^2} |\!|R_{\varphi}|\!|_{L^\infty} N^\eps(w).
\end{eqnarray}
Consequently, we can replace \eqref{estdtN1} in \eqref{estdtN0} and use the 
above estimates to get
\begin{eqnarray}
\label{estdtN3}
\frac{d}{dt} \int_{\R^d} \big( \mathcal{S}^\eps w, w \bigr) & = & 
\int_{\mathbb{R}^d} 4f'(|a^\eps|^2) (w, a^\eps) \Big( (w, \p_t a^\eps) - 
\big( u^\eps \cdot \nabla w + \frac{1}{2}\, w \, \nabla \cdot u^\eps, a^\eps \big) \Big) \\
 \nonumber & & + 2 \int_{\mathbb{R}^d} \partial_{t}\big[ f'(|a^\eps|^2)\big] (w, a^\eps)^2 + E_{1},
\end{eqnarray} 
where $E_{1}$ satisfies the estimate
\begin{eqnarray}
\label{dtNe1}
 E_{1} & \leq & C_M \bigg( 1 + 
\frac{1}{\e} \big|\!\big| R_{\varphi} \big|\!\big|_{W^{1,\ii}} + 
\frac{1}{\e^2} \big|\!\big| R_{\varphi} \big|\!\big|_{L^\ii} \bigg) N^\eps (w) \\
& & \nonumber - \frac{4}{\e} \int_{\R^d} f'(|a^\e|^2) (w,a^\e) (a^\e, i F^\e) + 
\int_{\mathbb{R}^d} (\eps \Delta w, i F^\eps).
\end{eqnarray}
To estimate the first integral in the right hand side of \eqref{estdtN3}, we use 
the equation \eqref{Ra} to get
\begin{align*}
4 \int_{\R^d} & f'(|a^\eps|^2) (w, a^\eps)\Big( (w, \partial_{t}a^\eps)- \big( u^\eps \cdot \nabla w + 
\frac{1}{2}\, w \, \nabla \cdot u^\eps, a^\eps \big) \Big) \\
& =4 \int_{\R^d} f'(|a^\eps|^2) (w, a^\eps) \Big( (w, R_{a}) - u^\eps \cdot \nabla (w, a^\eps) - 
(w, a^\eps) \nabla \cdot u^\eps \Big)\\
& = 4 \int_{\R^d} f'(|a^\eps|^2) (w, a^\e) (w, R_{a}) 
- 2 \int_{\R^d}f'(|a^\eps|^2)\, u^\eps \cdot \nabla \big( (w, a^\eps)^2 \big)
- 4 \int_{\R^d}f'(|a^\eps|^2) (w, a^\e)^2 \, \nabla \cdot u^\eps\\
& = 4 \int_{\R^d} f'(|a^\eps|^2) (w, a^\e) (w, R_{a}) 
+ 2 \int_{\R^d} (w, a^\eps)^2 \, u^\eps \cdot \nabla \big[ f'(|a^\eps|^2) \big]
- 2 \int_{\R^d} f'(|a^\eps|^2) (w, a^\e)^2\, \nabla \cdot u^\eps.
\end{align*}
To get  the last line, we have integrated by parts the second integral.
 Note that  the last term is  bounded by 
$ C_M N^\eps ( w )$, and, as for \iref{estdtN2} that the first integral is
 bounded by  
$\displaystyle{ \frac{C_M}{\e} |\!|R_{a} |\!|_{L^\infty} N^\eps ( w)}$. 
Consequently, we can replace the above identity in \eqref{estdtN3} to get
\beq
\label{estdtN4}
\frac{d}{dt} \int_{\R^d} \big( \mathcal{S}^\eps w, w \bigr) = 
\int_{\mathbb{R}^d} 2 (w, a^\eps)^2 \Big( \partial_{t} + u^\eps \cdot \nabla \Big) f'(|a^\eps|^2)  + E_{1}
+ E_{2} =: I + E_{1}+ E_{2},
\eeq
where $E_{2}$ is such that
\beq
\label{E2}
E_{2} \leq C_M \Big( 1 + \frac{1}{\e} \big|\!\big|R_{a} \big|\!\big|_{L^\infty} \Big) N^\eps ( w ).
\eeq
To estimate $I$, we use again the equation \eqref{Ra} which gives
$$\Big( \partial_{t} + u^\eps \cdot \nabla \Big) f'(|a^\eps|^2)
= 2f''(|a^\eps|^2) \big(a^\e, \partial_{t} a^\e+ u^\eps \cdot \nabla a^\e \big) = 
2 f''(|a^\eps|^2) \Big( R_{a} - \frac{1}{2 }\, a^\eps \, \nabla \cdot u^\eps , a^\eps \Big)$$
and hence we find
$$I \leq C \int_{\mathbb{R}^d} |a^\eps|^2\, \big| f''(|a^\eps|^2) \big|\, (w, a^\eps)^2 
+ 4\int_{\mathbb{R}^d} |a^\eps| \, |f''(|a^\eps|^2)|\, (w, a^\eps)^2 \, |R_{a}|.$$
To conclude, we shall use the assumption $(\BA)$. By defining  $n\in \N^*$ the first integer 
such that $f^{(n)}(0)\not = 0$, we see from Taylor expansion that
\be
\label{facto}
f'(\rho) = \rho^{n-1} q(\rho)
\ee
for some smooth positive function $q$ on $[0, + \infty)$. In particular, since $q > 0$, we have
$$\rho \mapsto \frac{ \rho f''(\rho)}{f'(\rho)} =n-1+ \rho \frac{q'(\rho)}{q(\rho)} 
\in \BC^\ii\big( [0,+\ii) \big),$$
which implies
\be
\label{regular}
\big| \rho f''(\rho) \big| \leq  C_M f'(\rho) \quad \quad \quad {\rm for} \ \ 0\leq \rho \leq M.
\ee
This yields 
$$ \int_{\mathbb{R}^d} |a^\eps |^2\, |f''(|a^\e|^2) | (w,a^\e)^2 
\leq C_M \int_{\mathbb{R}^d} (w,a^\e)^2 f'(|a^\e|^2) \leq C_M N^\eps (w),$$
where, again,  $C_M$ depends only on $M$. In a similar way, we also obtain
\begin{align*}
\int_{\mathbb{R}^d} (w,a^\e)^2 |a^\eps |\, | f''(|a^\e|^2) | \, |R_{a}| \leq 
& \ \big|\!\big| R_a \big|\!\big|_{L^\ii} 
\int_{\mathbb{R}^d} |w| \cdot \big| (w,a^\e) \big|\cdot |a^\eps |^2\, \big| f''(|a^\e|^2) \big|\\ 
\leq & \ \frac{C_M}{\e} \big|\!\big| R_a \big|\!\big|_{L^\ii} 
\int_{\mathbb{R}^d} \big( \e | w | \big) \bigg| (w,a^\e) \sqrt{f'(|a^\e|^2)} \bigg| \\
\leq & \ \frac{C_M}{\e} \big|\!\big| R_a \big|\!\big|_{L^\ii} N^\eps (w).
\end{align*}
Consequently, we have proven that
\beq
\label{estdtNI}
I \leq C_M\Big( 1 +\frac{1}{\e} \big|\!\big| R_a \big|\!\big|_{L^\ii} \Big) N^\eps(w).
\eeq
To get the result of Lemma \ref{dtN}, it remains to perform the $L^2$ estimate. Taking the 
$L^2$ scalar product of \iref{NLSwlin} with $iw$ and using that 
$$\displaystyle{(w, u^\e \cdot \nabla w + \frac{1}{2}\, w\, \nabla \cdot u^\e ) = 
\frac{1}{2}\, \nabla \cdot \big( |w|^2 u^\e \big)},$$
 we get
$$ \frac{ d}{dt} \Big( \frac{\eps^2 }{ 2} |\!|w|\!|_{L^2}^2 \Big) = 
\int_{\R^d} \e (F^\e, iw) + 2 \eps \int_{\R^d} f'(|a^\eps|^2) (w, a^\eps) ( a^\eps, iw).$$
Note that we have once again used that $R_{\varphi}$ is real-valued and hence that  
$(R_{\varphi }w, iw) = 0$. The first integral is clearly bounded by  $ N^\eps(w) + 
\big|\!\big| F^\e \big|\!\big|_{L^2}$
whereas for the second one, we have
$$  \int_{\mathbb{R}^d}  -2  \e f'(|a^\e|^2) (w,a^\e) (a^\e, iw)
 \leq  C_{M} \int_{\mathbb{R}^d}\Big(  f'(|a^\eps|^2) (w, a^\eps)^2 + \eps^2 |w|^2 \Big) 
 \leq C_{M} N^\eps(w(t)). $$
  As a consequence, we get
\beq
\label{dtNL2}
\frac{ d}{dt} \Big( \frac{\eps^2 }{ 2} \big|\!\big|w\big|\!\big|_{L^2}^2 \Big) 
\leq C_M N^\eps(w) + \big|\!\big|F^\eps \big|\!\big|_{L^2}^2.
\eeq
Finally, we can collect \eqref{dtNe1}, \eqref{estdtN4}, \eqref{E2}, \eqref{estdtNI} 
and \eqref{dtNL2} to get \eqref{lin0}. This completes the proof.\b

\subsection{Higher order estimates}

\ \indent Since our final aim is to prove Theorem \ref{Conver} by a fixed point argument, 
we also need to have $H^s$ estimates for $s$ sufficiently large for the solution of the 
linear equation \eqref{NLSwlin}. This is the aim of the following. Note that the term 
$- 2(w, a^\eps) f'(|a^\eps|^2) a^\eps$ in \eqref{NLSw} can be seen as a singular term 
with variable coefficients. Consequently, a crude way to get $H^s$ estimates is to apply 
$\eps^{|\alpha|} \partial^\alpha$ to the equation, the weight $ \eps^{| \alpha|}$ being 
used to compensate the singular commutator when we take the derivative of \eqref{NLSw}, 
and then to apply Lemma \ref{dtN} to the resulting equation. Nevertheless, it is possible 
to avoid the loss of $\eps^{|\alpha|}$ with more work by using more clever higher order 
modulated functionals. If $s\in \N$, $s\geq 2$, we define the following weighted norm, 
where $\alpha \in \N^d$ are multi-indices
\begin{eqnarray}
\label{nepsdef}
N^\eps_{s} (w) & \equiv &  
\sum_{ |\alpha | \leq s-1} N^\eps (\partial^\alpha w) + 
 K |\!| \mbox{Re } w|\!|_{H^{s-2}}^2  \\
\nonumber  & = &  { 1 \over 2} \eps^2 |\!| \nabla w |\!|_{H^{s-1}}^2  + 2 \sum_{ |\alpha | \leq s-1}
  \int f'(|a^\eps |^2) (\partial^\alpha w, a^\eps)^2   + K\big( \eps^2 |\!|w |\!|_{H^{s-1}}^2
   + |\!| \mbox{Re } w |\!|_{H^{s-2}}^2\big).
\end{eqnarray}
In this section, we shall use that 
$$ a^\eps = a^0 + \eps a^r$$
with $ a^0 $ real-valued  and 
$$ \sup_{ \eps \in (0, 1]} |\!| a^r |\!|_{L^\infty([0, T], W^{s, \infty})} \leq C. $$
Note that this allows to write
$$   \int_{\mathbb{R}^d} f'(|a^\eps |^2) (\partial^\alpha w, a^\eps)^2  \geq  { 1 \over 2 }
 \int_{\mathbb{R}^d} f'(|a^\eps |^2) (a^0)^2 |\mbox{Re } \partial^\alpha w|^2 - C
  \eps^2 |\!| \mbox{Re } 
 \partial^\alpha w |\!|_{L^2}^2
 $$
and  hence  by choosing $K$  sufficiently large ($K>C$) we get the lower bound
\beq
\label{equiv1}
N^\eps_s(w) \geq  \frac{1}{2} \sum_{|\alpha| \leq  s-1} N^\e\big( \p_x^\alpha w \big) 
 + \sum_{ |\alpha | \leq s-1 }
 \int_{\mathbb{R}^d} f'(|a^\eps |^2) (a^0)^2 |\mbox{Re } \partial^\alpha w|^2 \, dx 
\eeq
Note that we also have  the equivalence of norms:
\beq
\label{equivHs}
\big|\!\big| w \big|\!\big|_{H^s}^2 \leq \frac{2}{\eps^2} N_{s}^\eps(w), \quad
N^\eps_s(w) \leq C(|a^\eps|_{W^{s-1, \infty} }) \big|\!\big|w\big|\!\big|_{H^s}^2
 +  |\!| \mbox{Re } w |\!|_{H^{s-2}}^2.
\eeq

The main result of this section is :

\begin{theo}
\label{estiHS}
Let $0<T<\ii$, $s\in \N^*$, $f$ satisfying $(\BA)$ and $w\in \BC^1([0,T],H^s)$ 
 a solution of  \iref{NLSwlin} 
with $u^\e : [0,T] \times \R^d \to \R^d$ and 
$a^\e : [0,T] \times \R^d \to \C$ such that 
$$ M \equiv \sup_{0<\e<1} \bigg( \big|\!\big| u^\e \big|\!\big|_{L^\ii([0,T],W^{s+1,\ii}(\R^d))} + 
\big|\!\big| a^\e \big|\!\big|_{L^\ii([0,T],W^{s,\ii}(\R^d))} \bigg) < +\ii .$$
Assume finally that, for some $a^0 \in L^\ii([0,T],W^{s,\ii}(\R^d))$ {\rm real-valued}, 
$a^\e$ writes
\be
\label{aquasireel}
a^\e = a^0 + \BO_{W^{s,\ii}}(\e)
\ee
uniformly on $[0,T]$. Then, there exists $C$, depending only on $d$, $f$ and $M$, such that
$$ \frac{d}{dt} N^\eps_s \big( w(t) \big) \leq 
C \bigg( 1+ \frac{1}{\eps} \big|\!\big| R_a (t) \big|\!\big|_{L^{\ii}} + 
\frac{1}{\eps^2} \big|\!\big| R_\varphi (t) \big|\!\big|_{W^{s-1,\ii}}
 \bigg) N^\eps_s \big( w(t) \big) 
+ C \big|\!\big| F^\e(t) \big|\!\big|_{H^s}^2 
+ \frac{C}{\e^2} \big|\!\big| {\rm Im}\, F^\e(t) \big|\!\big|_{H^{s-1}}^2.$$
\end{theo}

\begin{rem}\rm
In view of \iref{aquasireel}, $a^\e$ is real up to $\BO(\e)$, hence, in the integral 
in the right-hand side of \iref{lin0}, the real and imaginary parts of $F^\e$ do 
not play the same role. This explains that the estimate is better for ${\rm Re}\, F^\e$ 
than for ${\rm Im}\, F^\e$. As a matter of fact, for $s=1$, Theorem \ref{estiHS} follows 
immediately from Lemma \ref{dtN} and \iref{aquasireel}.
\end{rem}
\subsection{Proof of Theorem \ref{estiHS}}

We estimate separately the two terms  in $N^\eps_s(w)$, when $s\geq 2$ (otherwise, the
 result  follows 
from Lemma \ref{dtN} as we have seen). Let us set 
$$ \Sigma(w) \equiv 
 |\!|\mbox{Re } w |\!|_{H^{s-2}}^2.$$
 Note that  we have
 \beq
 \label{sigmaequiv} \Sigma(w) \leq N^\eps_{s}(w).
 \eeq
In the proof, $C$ is a constant depending only on $d$, $f$ and $M$.

We shall first 
 prove that 
\be
\label{deuz}
\frac{d}{dt} \Sigma(w) \leq 
C \Big( 1 + 
 \frac{1}{\e^2} \big|\!\big| R_\varphi \big|\!\big|_{W^{s-2,\ii}} \Big) N^\eps_s\big( w \big) 
+ C \big|\!\big| F^\e \big|\!\big|_{H^{s-2}}^2 + 
\frac{C}{\e^2} \big|\!\big| {\rm Im}\, F^\e \big|\!\big|_{H^{s-2}}^2.
\ee
\ \\
For $\alpha \in \N^d$, we have
\begin{eqnarray}
\label{dwa}
\partial_{t} \big( \partial^\alpha  w \big) + u^\e \cdot \nabla \big( \p^\alpha w \big) 
= & & \frac{i \e}{2} \Delta (\p^\alpha w) - \frac{i}{\e} \p^\alpha F^\e - \frac{i}{\eps} \p^\alpha
\big(R_{\varphi } w \big)  \\ 
\nonumber & &
- \frac{2i}{\eps} \p^\alpha\big( f'(|a^\eps|^2) (a^\eps, w)  a^\eps\big) 
- \big[ \p^\alpha, u^\e \cdot \nabla \big] w - 
\frac{1}{2}\, \p^\alpha \big( w \nabla \cdot u^\e \big).
\end{eqnarray}
Next,  by taking the real part of \eqref{dwa},  we get
\begin{eqnarray*}
\partial_{t} \big( \partial^\alpha \mbox{Re }  w \big) + u^\e \cdot \nabla \big( \p^\alpha 
\mbox{Re }w \big) 
 & =  & 
- \big[ \p^\alpha, u^\e \cdot \nabla \big] \mbox{Re  } w - 
\frac{1}{2}\, \p^\alpha \big( \mbox{Re } w \,  \nabla \cdot u^\e \big) + \mathcal{R}^\eps
\end{eqnarray*}
where
\beq
\label{realeq} \mathcal{R}^\eps = 
\mbox{Re}\Big( \frac{i \e}{2} \Delta (\p^\alpha w) - \frac{i}{\e} \p^\alpha F^\e - \frac{i}{\eps} \p^\alpha
\big(R_{\varphi } w \big)  
- \frac{2i }{\eps} \p^\alpha\big( f'(|a^\eps|^2) (a^\eps, w)  a^\eps \big) \Big).
\eeq 
  By using  \eqref{aquasireel}, we  have 
    $$\mbox{Im }  \partial^\gamma a^\eps = \mathcal{O}(\eps), \quad \forall \gamma, \, 
     |\gamma | \leq |\alpha| $$
    and
  \beq
  \label{areel2}
  |(\partial^\beta a^\eps, \partial^\gamma w)| \leq C |\mbox{Re }\partial^\gamma w | + \eps |\partial^\gamma w|
  \eeq
  for every $\beta,$ $\gamma$.
  Consequently, we immediately obtain for every $\alpha$, $|\alpha |\leq s-2$,    
 \begin{eqnarray*}
  |\!| \mathcal{R}^\eps |\!|_{L^2} 
  & \leq &  C \Big(  \eps |\!| w |\!|_{H^s}
   + {  |\!| R_{\varphi} |\!|_{W^{s-1, \infty}}  \over \eps^2 } |\!|w|\!|_{H^{s-2}}
    +  |\!| \mbox{Re } w |\!|_{H^{s-2}} + \eps  |\!|   w |\!|_{H^{s-2}}
   \Big) +  { 1 \over \eps} |\!|\mbox{Im }F^\eps|\!|_{H^{s-2}}\\
    & \leq & C\big( 1 +  {  |\!| R_{\varphi} |\!|_{W^{s-2, \infty}}  \over \eps^2 } ) N^\eps_{s}(w)^{1 \over 2} +   { 1 \over \eps} |\!|\mbox{Im }F^\eps|\!|_{H^{s-2}}.
    \end{eqnarray*}
    Consequently, the standard $L^2$ energy estimate for \eqref{realeq} gives
\beq
\label{Realw}
 { d \over dt } |\!| \mbox{Re }\partial^\alpha |\!|_{L^2}^2 \leq   C\big( 1 +  {  |\!| R_{\varphi} |\!|_{W^{s-1, \infty}}  \over \eps^2 } ) N^\eps_{s}(w) +   { 1 \over \eps^2} |\!|\mbox{Im }F^\eps|\!|_{H^{s-2}}^2.\eeq
Note that  we have used that 
    $$ \int   \Big( u^\e \cdot \nabla \big( \p^\alpha 
\mbox{Re }w \big), \partial^\alpha \mbox{Re } w \Big)
 =- {Ê 1 \over 2 } \int( \nabla \cdot u^\eps ) |\partial^\alpha \mbox{Re } w |^2.$$ 
     Consequently, \eqref{deuz} is proven.
  \bigskip 
   
   The next step is to estimate $N^\eps (\partial^\alpha w )$ for $ |\alpha | \leq s-1$.
    By applying $\p^\alpha$ to \iref{NLSwlin}, we get 
\be
\label{eqrase}
i\e \frac{\p (\p^\alpha w)}{\p t} + \BL^\e\big( \p^\alpha w \big) 
= R_{\varphi} \partial^\alpha w +\tilde{F}^\e,
\ee
where
$$ \tilde{F}^\eps \equiv \mathcal{C}^\alpha + \BD^\alpha + \p^\alpha F^\eps
+ [\p^\alpha , R_\varphi] w,$$
with
\begin{eqnarray*}
\mathcal{C}^\alpha & \equiv & 2 \,\p^\alpha\Big( f'(|a^\eps |^2) a^\eps (w, a^\eps) \Big) 
- 2 f'(|a^\eps|^2)(\p^\alpha w, a^\eps)a^\eps, \\
\mathcal{D}^\alpha & \equiv & 
 -i \e \big[ \p^\alpha, u^\e \cdot \nabla \big] w 
-\frac{i \e}{2} \big[ \p^\alpha,\, \nabla \cdot u^\e \big] w.
\end{eqnarray*}
To estimate $N^\eps(\partial^\alpha w)$, we shall use Lemma \ref{dtN}.
Towards this, we need  to estimate the commutators in the right hand side
 of \eqref{eqrase}.
For 
$|\alpha| \leq s-1$, the following estimates hold for $\BC^\alpha$ and $\BD^\alpha$:
\begin{align}
\label{estC3}
\big|\!\big| [\partial^\alpha, R_{\varphi} ] w \big|\! \big|_{H^1 }^2 & \ \leq 
C \big|\!\big|R_{\varphi} \big|\!\big|_{W^{s,\infty} }^2\big|\! \big| w \big|\! \big|_{H^{s}}^2 
\leq \frac{C}{\e^2} \big|\!\big|R_{\varphi} \big|\!\big|_{W^{s,\infty} }^2 N^\eps_{s}(w),\\
\label{estC1}
\big|\!\big|\BD^\alpha \big|\!\big|_{H^1}^2 & \ \leq 
C\, \eps^2 \big|\! \big| w \big|\! \big|_{H^{s}}^2 \leq C N^\eps_{s}(w),\\
\label{Dalf}
\big|\!\big| \big( i f'(|a^\eps |^2)^{ 1 \over 2 }  a^\e, \BD^\alpha \big) \big|\!\big|_{L^2}^2 & \ \leq C \e^2 N^\e_s(w),\\
\label{CalfH1}
\big|\!\big| \BC^\alpha \big|\!\big|_{H^1}^2 & \ \leq C  N^\eps_{s}(w), \quad \quad \quad\\
\label{Calfscalia}
\big|\!\big| (ia^\e,\BC^\alpha) \big|\!\big|_{L^2}^2 & \ \leq C \e^2 N^\eps_{s}(w).
\end{align}
The estimates \iref{estC3} and \iref{estC1} follow easily from \eqref{equivHs}. For \iref{Dalf}, 
we note that
\begin{align*}
\frac{1}{\e} \big( i a^\e, \BD^\alpha \big) = & \ 
- \big( a^\e, [\p^\alpha , u^\e \cdot \nabla ] w \big) 
- \frac{1}{2} \big( a^\e, [\p^\alpha , \nabla \cdot u^\e ] w \big)\\
= & \ - \sum_{\gamma < \alpha} \left(\begin{array}{c} \alpha\\ \gamma \end{array}\right) 
\big( \p^{\alpha - \gamma} u^\e \big) \cdot \big( a^\e, \nabla \p^\gamma w \big) - 
\frac{1}{2}\ \sum_{\gamma < \alpha} \left(\begin{array}{c} \alpha\\ \gamma \end{array}\right) 
\p^{\alpha - \gamma} \big( \nabla \cdot u^\e\big) \big( a^\e,\p^\gamma w \big)
\end{align*}
 since  $u^\eps$ is real.   Next, we can use   \eqref{aquasireel} 
  and \eqref{areel2} again. 
 In particular,  in the above
 expansion,  the terms 
$\big( a^\e,\p^\gamma w \big)$  are bounded in 
$L^2$ by $\Sigma(w) + \eps^2 |\!| w|\!|_{H^{s-2}}^2$  and thus by $N^\eps_{s}(w)$.
Similarly, 
  the terms  $\big( a^\e, \nabla \p^\gamma w \big)$ are
  bounded in $L^2$ by  $N^\eps_{s}(w)$ if $|\gamma |\leq  s-3.$
Consequently, we get
$$\big|\!\big| \big( i f'(|a^\eps |^2)^{ 1 \over 2 }  a^\e, \BD^\alpha \big) \big|\!\big|_{L^2}^2
 \leq C \Big( \sum_{|\beta|=s-1}
  \int_{\mathbb{R}^d}  f'(|a^\eps |^2)( \partial^\beta w, a^\eps)^2 +
    N^\eps_{s}(w) \Big) \leq C N^\eps_{s}(w)$$ 
    which yields \iref{Dalf}.
Next, we  turn to $\BC^\alpha$.  The   Leibnitz formula  gives
\be
\label{develoC}
\BC^\alpha = 
\sum_{\scriptsize{\begin{array}{c} \t{\alpha}<\alpha,\\  \t{\alpha}  + \beta + \lambda +\mu =\alpha \end{array}}} 
* \ \p^\lambda \big[ f'(|a^\e|^2) \big] \big( \p^{\t{\alpha}} w , \p^\beta a^\e \big) \p^\mu a^\e,
\ee
where $*$ is a  real coefficient depending only on $\t{\alpha}$, $\beta$, $\lambda$ and $\mu$.
Since 
$|\t{\alpha}| \leq |\alpha| -1 \leq s-2$, we can use again \eqref{aquasireel}
 through \eqref{areel2} to get
 that
 $$ |\!| \mathcal{C}^\alpha |\!|_{L^2}^2 \leq C\Big( \Sigma(w) + \eps^2 |\!|w |\!|_{H^s}^2
  \Big) \leq C N^\eps_{s}(w).$$
Since  $(ia^\eps, \partial^\mu a^\eps)= \mathcal{O}(\eps)$ thanks to \eqref{aquasireel},
 we  also get \eqref{Calfscalia}.  
For the $H^1$ norm,  the same argument yields
$$ |\!| \mathcal{C}^\alpha |\!|_{H^1}^2 \leq C\Big( \Sigma(w) + \eps^2 |\!|w |\!|_{H^s}^2
  + \sum_{\scriptsize{\begin{array}{c} {|\gamma|}=s-1,\\  | \beta + \lambda +\mu| =1 \end{array}}} \int_{\mathbb{R}^d}  \big| \p^\lambda \big[ f'(|a^\e|^2) \big] \big( \p^{\gamma} w , \p^\beta a^\e \big) \p^\mu a^\e  \big|^2
     \Big).$$
  To estimate the last sum, we first consider the terms   with  $\beta =0$.
   They are always bounded by 
  $$   C\int_{\mathbb{R}^d}   \big[ f'(|a^\e|^2) + |a^\eps|^2 |f''(|a^\eps |^2)| \big] \big( \p^{\gamma} w , a^\e \big)^2$$
  with $|\gamma| = s-1$ and hence, thanks to \eqref{regular},
   they are bounded by 
  $$ C\int_{\mathbb{R}^d}    f'(|a^\e|^2) \big( \p^{\gamma} w , a^\e \big)^2$$
   and hence by $N^\eps_{s}(w)$.
Next, we consider the terms with $|\beta|=1$.  Since
 $\lambda = \mu = 0$,  we have
 to estimate terms like
$$ \mathcal{T}_{1}=  \int_{\mathbb{R}^d}   f'(|a^\e|^2) \big( \p^{\gamma} w ,\partial^\beta a^\e \big)^2 |a^\eps|^2.$$
By using again \eqref{aquasireel} and \eqref{areel2}, we get
$$ \mathcal{T}_{1} \leq C \int_{\mathbb{R}^d}   f'(|a^\e|^2) |a^0|^2 |\mbox{Re } \p^{\gamma} w|^2
  + C \eps^2 |\!| w|\!|_{H^{s-1}}^2$$
  and hence, by using \eqref{equiv1}, we finally obain
 $$  \mathcal{T}_{1} \leq C N^\eps_{s}(w).$$
Consequently, \eqref{CalfH1} is proven. This ends the estimates
 of the commutators.

We are now  able to  establish:
\begin{align}
\label{cestmieux}
\frac{d}{dt} N^\e \big( \p^\alpha w \big) \leq C &\ 
\Big( 1+ \frac{1}{\e^2} \big| \!\big| R_{\varphi} \big|\!\big|_{W^{s-1,\ii}} + 
\frac{1}{\e} \big|\!\big| R_a \big|\!\big|_{L^\ii} 
\Big) N^\eps_{s}(w)\nonumber \\
& \ + \big| \!\big| F^\eps \big|\!\big|_{H^s}^2 + 
\frac{C}{\e^2} \big|\!\big| {\rm Im}\, F^\e \big|\!\big|_{H^{s-1}}^2.
\end{align}

Indeed, from Lemma \ref{dtN}, we deduce
\begin{align}
\nonumber \frac{d}{dt} N^\e \big( \p^\alpha w \big) \leq 
C &\ \Big( 1+ \frac{1}{\e} \big| \!\big| R_{\varphi} \big| \!  \big|_{W^{1,\ii}} + 
\frac{1}{\e} \big|\!\big| R_a \big|\!\big|_{L^\ii} + 
\frac{1}{\e^2} \big|\!\big| R_\varphi \big|\!\big|_{L^\ii} \Big) N^\e \big( \p^\alpha w \big)\\
\label{term1}& \ + \big| \!\big| \tilde{F}^\eps \big|\!\big|_{L^2}^2 
+ \frac{4}{\e} \int_{\mathbb{R}^d } f'(|a^\e|^2) (\p^\alpha w,a^\e) \big( i a^\e, \tilde{F}^\e \big) 
- \int_{\mathbb{R}^d} (i \eps \Delta \p^\alpha w, \tilde{F}^\eps ),
\end{align}
To estimate the right-hand side of \eqref{term1}, we first estimate 
$ \big| \!\big| \tilde{F}^\eps \big|\!\big|_{L^2}^2$.
Combining \iref{estC3} and \iref{estC1} with \iref{CalfH1}, we infer
\beq
\label{FL2}
\big| \!\big| \tilde{F}^\eps \big|\!\big|_{L^2}^2 
\leq \big| \!\big| F^\eps \big| \!\big|_{H^{s-1}}^2 + 
C\Big( 1 + \frac{1}{\e^2} \big|\!\big|R_{\varphi}\big |\!\big|_{W^{s-1, \infty}}^2 \Big) N^\eps_{s}(w).
\eeq
Next, we  turn to the term
$$\frac{4}{\e} \int_{\mathbb{R}^d } f'(|a^\e|^2) (\p^\alpha w,a^\e) \big( i a^\e, \tilde{F}^\e \big) = 
\frac{4}{\e} \int_{\mathbb{R}^d} f'(|a^\e|^2) (\p^\alpha w,a^\e) \big( i a^\e, 
\mathcal{C}^\alpha + \BD^\alpha + \p^\alpha F^\eps + [\p^\alpha , R_\varphi] w \big),$$
which splits as four integrals. For the first one, by \iref{Calfscalia} and Cauchy-Schwarz:
$$ \frac{4}{\e} \int_{\mathbb{R}^d} f'(|a^\e|^2) (\p^\alpha w,a^\e) \big( i a^\e, \mathcal{C}^\alpha \big) 
\leq C \Big( \int_{\mathbb{R}^d} f'(|a^\eps|^2) \big(\partial^\alpha w,
 a^\eps \big)^2 \Big)^{1 \over 2 }\, N^\e_s(w)^{\frac{1}{2}}  \leq CN^\e_s(w).$$
For the second one, we use \iref{Dalf} and Cauchy-Schwarz, which gives
$$ \frac{4}{\e} \int_{\mathbb{R}^d} f'(|a^\e|^2)^{ 1 \over 2} (\p^\alpha w,a^\e) \big( i
f'(|a^\e|^2)^{ 1 \over 2} a^\e, \BD^\alpha \big) 
\leq CN^\e_s(w).$$
For the third integral, we simply write, using once again \iref{aquasireel}
$$ \frac{1}{\e} \big| \!\big| \big( i a^\e,\p^\alpha F^\eps \big) \big| \!\big|_{L^2} 
\leq C \big|\!\big|  F^\e \big|\!\big|_{H^{s-1}} 
+ \frac{C}{\e} \big|\!\big| {\rm Im}\, F^\e \big|\!\big|_{H^{s-1}},$$
which yields by Cauchy-Schwarz
$$ \frac{4}{\e} \int_{\mathbb{R}^d} f'(|a^\e|^2) (\p^\alpha w,a^\e) \big( i a^\e,\p^\alpha F^\eps \big) 
\leq C N^\eps_{s}(w) + C\big|\!\big|  F^\e \big|\!\big|_{H^{s-1}}^2 
+ \frac{C}{\e^2} \big|\!\big| {\rm Im}\, F^\e \big|\!\big|_{H^{s-1}}^2.$$
Finally, for the fourth integral, we have by \iref{estC3}
$$ \frac{4}{\e} \int_{\mathbb{R}^d} f'(|a^\e|^2) (\p^\alpha w,a^\e) \big( i a^\e, [\p^\alpha ,R_\vp] w \big) 
\leq \frac{C}{\e} \big|\!\big| R_\vp \big|\!\big|_{W^{s-1,\ii}} N^\e(w).$$
 By summing these estimates, we find 
\be
\label{Terme2}
\frac{4}{\e} \int_{\mathbb{R}^d } f'(|a^\e|^2) (\p^\alpha w,a^\e) \big( i a^\e, \tilde{F}^\e \big) 
\leq C\Big( 1 + \frac{1}{\eps} \big|\!\big|R_{\varphi} \big|\!\big|_{W^{s-1, \ii}} \Big) N^\eps_{s}(w)
+ C\big|\!\big| F^\e \big|\!\big|_{H^{s-1}}^2 
+ \frac{C}{\e^2} \big|\!\big| {\rm Im}\, F^\e \big|\!\big|_{H^{s-1}}^2.
\ee
Finally, we  handle the term
$$ - \int_{\mathbb{R}^d} \big(i \eps \Delta \p^\alpha w, \tilde{F}^\eps \big)= 
- \int_{\mathbb{R}^d} \big(i \eps \Delta \p^\alpha w, 
\BC^\alpha + \BD^\alpha + \p^\alpha F^\eps + [\p^\alpha , R_\varphi] w\big).$$
By using an integration by parts, we  have 
\begin{align*}
- \int_{\mathbb{R}^d} \big(i \eps \Delta \p^\alpha w, 
 \tilde{F}^\eps) \leq & \ 
\big|\!\big|\mathcal{C}^{\alpha} \big|\!\big|_{H^1} ^2 + 
\big|\!\big|\mathcal{D}^{\alpha} \big|\!\big|_{H^1} ^2
+ \big|\!\big| [\partial^\alpha, R_{\varphi} ] w \big|\!\big|_{H^1}^2 
+ \big|\!\big|F ^\eps \big|\!\big|_{H^s}^2 + C N^\eps_{s}(w)\\
 \leq & \ \big|\!\big|F ^\eps \big|\!\big|_{H^s}^2 
+ C\Big( 1 + \frac{1}{\eps^2} \big|\!\big|R_{\varphi} \big|\!\big|_{W^{s-1, \infty}} \Big) N^\eps_{s}(w)
\end{align*}
 thanks to  \eqref{estC1},  \eqref{estC3} and \eqref{CalfH1}. 
 
 Consequently, we can collect  the last estimate and  \iref{FL2}, \iref{Terme2},  \iref{term1} 
  to get  \eqref{cestmieux}. This ends the proof of Theorem \ref{estiHS}.

\section{Construction of WKB expansions}
\label{sectionWKB}

\ \indent In this section, we construct an approximate solution of \iref{NLS} using a WKB 
expansion. The first step is to prove the local existence of smooth solutions of the limit 
hydrodynamical system.

\subsection{Well-posedness of the limit system}

\ \indent We consider the system
\be
\label{hydro}
\left\{\begin{array}{ll}
\displaystyle{ \partial_{t} a + u \cdot \nabla a + \frac{1}{2}\, a \, \nabla \cdot u}  = 0\\ \\
\displaystyle{\partial_{t} u + u \cdot \nabla u + \nabla \big( f(a^2) \big) }  = 0,
\end{array}\right.
\ee
which is only weakly hyperbolic, with the pressure law $f$ satisfying assumption $(\mathcal{A})$, 
and the initial condition $(a,u)_{|t=0}=(a^0,u^0)$.

\begin{theo}
\label{theoeulernh}
Assume that $f$ satisfies assumption $(\mathcal{A})$, and denote by $n\in \N^*$ 
the first integer such that $f^{(n)}(0) \neq 0$. Then, for $s> 2 +d/2$, for every 
$(a_{0}, u_{0}) \in H^{s-1} \times H^s$ with $a_{0} \in \mathbb{R} $ and $a_0^n \in H^s$, there 
exists $T>0$ and a unique solution $(a,u)$ of \eqref{hydro} such that 
$(a, a^n, u) \in \mathcal{C}([0, T ], H^{s-1} \times H^s \times H^s) \cap 
\mathcal{C}^1([0, T ], H^{s-2} \times H^{s-1} \times H^{s-1})$.
\end{theo}

Let us remark that if $n=1$, then $f'(0)>0$ and thus $f'>0$ in $[0,+\ii)$ (by $(\mathcal{A})$). 
In this case, \iref{hydro} is symmetrizable (using the symmetrizer $S=$ 
diag$\big(1, \frac{1}{4}f'(a^2), ... , \frac{1}{4}f'(a^2) \big)$ used in $\cite{G}$) and 
the existence and uniqueness for \iref{hydro} follows easily.

\subsection*{Proof of Theorem \ref{theoeulernh}}
The first step is to rewrite the system by using more convenient unknowns. At first, 
we notice that thanks to $(\mathcal{A})$, we can write $f$ under the form
$$ f(\rho)= \rho^n \tilde{f}(\rho),$$
with $\tilde{f}$ smooth on  $[0,+\ii)$ and such that $\tilde{f}(0)\neq 0$. 
Next, since we have by assumption $f(0)=0$ and $f'(\rho)>0$ for $\rho \neq 0$, 
we also have that $f(\rho)>0$ for $\rho>0$. This implies that $\tilde{f}(\rho)>0$ 
for $\rho\geq 0$. This allows to define a  smooth function  $h$
  on  $\mathbb{R}$ by 
\be
\label{hdef}
{\rm h}(a) \equiv a \big[ \tilde{f}(a^2) \big]^{\frac{1}{2 n}}.
\ee
Note that   h$(a) \neq 0$ for $a\neq 0$.
It is usefull to notice that  we can also write h under the form
$$ {\rm h}(a)=  \mbox{sgn}(a)\, f(a^2)^{ 1 \over 2 n }$$
 and hence that we have
$$
{\rm h}(a)^{2n} =  f(a^2).
$$
 Furthermore, since $f'>0$ and $\tilde{f}(0)>0$, 
we deduce that h$'(a)>0$  for $a \neq 0$  and  that h$'(0) = \big[ \tilde{f}(0) \big]^{\frac{1}{2 n}}>0$, 
so that h$'>0$ on $\mathbb{R}$. Thus h is a smooth diffeomorphism from $\mathbb{R}$
 to 
h$\big(\mathbb{R} \big)$. In particular, this allows to define a smooth positive function $c$ 
on  h$\big( \mathbb{R} \big)$ such that
$$ \frac{1}{2}\, a {\rm h}'(a) = {\rm h}(a) \, c \big( {\rm h}(a) \big), \quad \forall a \in
\mathbb{R}.$$
With this definition, $(h,u)$, with $h \equiv$ h$(a)$, solves the system
\be
\label{hydroh}
\left\{\begin{array}{ll}
\displaystyle{ \partial_{t} h + u \cdot \nabla h + h c( h ) \, \nabla \cdot u}  = 0\\ \\
\displaystyle{\partial_{t} u + u \cdot \nabla u + \nabla \big( h^{ 2n} \big)}  =0.
\end{array}\right.
\ee
Since $a$ is in $H^s$ if and only if $h$ is in $H^s$, we shall prove local existence of 
smooth solution for the weakly hyperbolic system \eqref{hydroh}. 
As we shall see below, the nonlinear symmetrization method of \cite{MUK}
 does not allow to reduce \eqref{hydroh} to a symmetric or
  symmetrizable system with smooth coefficients except
   in the case that $c(h)= \tilde{c}( h^n)$ with $\tilde{c}$ smooth . Nevertheless, it 
    will be  still
   possible to use the same idea  to prove the existence of an energy estimate
    with loss
    for the system \eqref{hydroh}. When we are in such a situation, 
the simplest way to 
construct a solution is to use the vanishing viscosity method. Indeed, this approximation 
method allows to preserve the nonlinear energy estimate verified by \eqref{hydroh}. We 
 thus consider for $\ep>0$ the system
\be
\label{hydrov}
\left\{\begin{array}{ll}
\displaystyle{\partial_{t} h_\ep + u_\ep \cdot \nabla h_\ep + h_\ep c(h_\ep) \nabla \cdot u_\ep}  = 
\ep\, \Delta h_\ep\\ \\
\displaystyle{\partial_{t} u_\ep + u_\ep \cdot \nabla u_\ep + \nabla \big( h_\ep^{2n} \big)}  
= \ep\, \Delta u_\ep.
\end{array}\right.
\ee
The local existence of smooth solutions for this parabolic system is very easy to obtain. 
Moreover, we note that $h_\ep$ remains nonnegative if the initial datum $h_{|t=0}$ is nonnegative. 
In the following, we shall only prove an $H^s$ energy estimate independent of $\ep$ for this system 
which ensures that the solution  remains smooth on an interval of time independent of $\ep$. 
The final step which consists in using the uniform bounds to pass to the limit when $\ep$ goes to zero to 
get a solution of \eqref{hydroh} is very classical and hence will not be detailled. 
In the proof of the energy estimates, we shall omit the subscript $\ep$ for notational 
convenience.

As in the work of \cite{MUK}, we introduce the unknown $H \equiv h^n = a^n \tilde{f}(a^2)^{\frac{1}{2}}$. 
Note that by definition of $h$, $H$ is in $H^s$ as soon as $a^n$ is in $H^s$. We get for $(H, u)$ 
the system
\be
\label{systHu}
\left\{\begin{array}{ll}
\displaystyle{\p_t H + u \cdot \nabla H + n H c(h)\, \nabla \cdot u}  = \ep\, n  h^{n-1} \Delta h = 
\ep\, \Big( \Delta H - n(n-1) h^{n-2} |\nabla h |^2\Big)\\ \\
\displaystyle{\p_t u \ + u \cdot \nabla u \ + \ 2 H \nabla H}  = \ep\, \Delta u.
\end{array}\right.
\ee
Note that it does not seem possible to get a classical hyperbolic symmetric system 
(in the case $\ep =0$) involving only $H$ and $u$ as in the case of  homogeneous 
pressure laws considered in \cite{MUK}. Indeed, the coefficient $c(h)= c(H^{\frac{1}{n}})$ 
is not (in general) a smooth function of $H$. Nevertheless, it will be possible to prove 
that the system with unknowns $(h, H, u)$ though only weakly hyperbolic (when $\ep = 0$) 
satisfies an energy estimate. We notice that the symmetrizer
$$ S \equiv \mbox{diag}\Big( 1, \frac{n}{2} c(h) I_{d} \Big),$$
which is positive since $c(h)$ is positive, symmetrizes the first order part of 
\eqref{systHu}. We shall first perform an $H^s$ energy estimate ($s>2 + d/2$) on 
\eqref{systHu} but we have to track carefully the dependence on $h$ in the energy estimates.

To prove our $H^s$ energy estimate, we shall make an extensive use of the following 
classical (see \cite{Taylor} for example) tame estimates
\begin{eqnarray}
\label{sob1} & & 
\big|\! \big| f g \big|\! \big|_{H^k} \leq 
C_{k} \Big( \big|\! \big|f\big|\! \big| _{L^\infty} \big|\! \big| g \big|\! \big|_{H^k} + 
\big|\! \big|f \big|\! \big|_{H^k} \big|\! \big|g \big|\! \big|_{L^\infty} \Big),\\
\label{sob2}& & 
\big|\! \big| \p^\alpha (f g ) - f \partial^\alpha g \big|\! \big|_{L^2} \leq C_k
\Big( \big|\! \big| f \big|\! \big|_{H^k} \big|\! \big| g \big|\! \big|_{L^\infty} + 
\big|\! \big| \nabla f \big|\! \big|_{L^\infty} \big|\! \big|g \big|\! \big|_{H^{k-1}} \Big), 
\quad |\alpha | \leq k,\\
\label{sob3}& &
\big|\! \big|F(u)\big|\! \big|_{H^k} \leq C( \big|\! \big|u \big|\! \big|_{L^\infty}) 
( 1 + \big|\! \big|u\big|\! \big|_{H^k})
\end{eqnarray}
 if $F$ is smooth and such that $F(0)=0$.
 
    At first, we notice that $(\p^\alpha H, \p^\alpha u)$ for $|\alpha| \leq s$ 
solves the system
$$\left\{\begin{array}{ll}
\displaystyle{\p_{t} \p^\alpha H+ u \cdot \nabla \p^\alpha H + n c(h) \big( \nabla \cdot u \big) \p^\alpha H } 
& = \ep \, \Big( \Delta \p^\alpha H -n(n-1) \p^\alpha( h^{n-2} |\nabla h |^2) \Big)\\ 
 & \ \ \ - [\p^\alpha, u ] \cdot \nabla H - n [\p^\alpha,  Hc(h)] \nabla \cdot u\\ \\
\displaystyle{\partial_{t} \p^\alpha u \ \, + u \cdot \nabla \p^\alpha u \ + \ \ \ 2 H \nabla \p^\alpha H} & = 
\ep \, \Delta \p^\alpha u - [\p^\alpha, u ] \cdot \nabla u - [\p^\alpha, 2 H] \nabla H.
\end{array}\right.$$
By using \eqref{sob2} to estimate in $L^2$ the commutators 
in the right hand-side,
we get in a classical way 
by integration by parts
\begin{eqnarray}
\label{eHs1}& & \frac{d}{dt}\Big[ \frac{1}{2}\int_{\R^d} |\p^\alpha H|^2 + 
\frac{n}{2} c(h) |\p^\alpha u |^2 \Big] + 
\ep\, \int_{\R^d} |\nabla \p^\alpha  H|^2 + \frac{n}{2} c(h) |\nabla \p^\alpha u |^2 \\
\nonumber & & \leq 
C_0 \Big( \big| \! \big| (h,u) \big| \! \big|_{W^{ 1, \ii}} \Big) 
|\!|V|\!|_{H^s}^2 + \mathcal{C}^\alpha + \ep\, \mathcal{D}^\alpha + \mathcal{R}^\alpha,
\end{eqnarray}
where $V \equiv (H, u)$, $C_0$ is a non-decreasing function depending only on $f$, $s$ and $d$, and
\begin{eqnarray*}
& & \mathcal{C}^{\alpha} \equiv - n \int_{\R^d} (\p^\alpha  H )\, [\p^\alpha, H c(h) ] (\nabla \cdot u),\\
& & \mathcal{D}^{\alpha} \equiv - \frac{n}{2} \int_{\R^d} 
c'(h) \big( (\nabla h \cdot \nabla ) \p^\alpha u \big) \cdot \p^\alpha u 
- n(n-1) \int_{\R^d} \p^\alpha \big( h^{n-2} |\nabla h |^2\big) \, \p^\alpha H, \\
& & \mathcal{R}^\alpha \equiv \frac{n}{4} \int_{\R^d} c'(h) \p_t h |\p^\alpha u |^2.
\end{eqnarray*}
We have singled out the three terms above since they are the ones involving $h$ which must
be estimated with care. Note that the estimate of $\mathcal{C}^\alpha$ will be crucial since 
this term involves high order derivatives of $h$. Next, we can integrate \eqref{eHs1} in time, 
sum the estimates for $|\alpha| \leq s$ and use that $c(h)>0$, hence 
$n c(h) /2 \geq \frac{1}{C_1(|\!| h|\!|_{L^\ii})}$ to obtain
\begin{eqnarray}
\label{eHs2}
& & |\!| V(t)|\!|_{H^s}^2 + \ep \int_{0}^t |\!|\nabla V(\tau) |\!|_{H^s}^2 \, d\tau \\
& &\nonumber \leq C_1 \big( |\!|h|\!|_{L^\infty} \big) \Big( |\!| V(0)|\!|_{H^s}^2 
+ \int_{0}^t C_0 \big( |\!|(h,u)(\tau)|\!|_{W^{1,\ii}} \big) |\!|V(\tau)|\!|_{H^s}^2 + \mathcal{C}(\tau) + \ep \mathcal{D}(\tau) + \mathcal{R}(\tau) 
\, d\tau \Big),
\end{eqnarray}
with
$$ \mathcal{C} \equiv \sum_{ |\alpha | \leq s} \mathcal{C}^\alpha, \quad 
\mathcal{D} \equiv \sum_{ |\alpha | \leq s} \mathcal{D}^\alpha,\quad  
\mathcal{R} \equiv \sum_{ |\alpha | \leq s} \mathcal{R}^\alpha.$$

\noindent {\bf Estimate for $\boldsymbol{\BC}$.} We claim that
\beq
\label{eCalpha}
\mathcal{C} \leq 
C_0( |\!|(h,u)|\!|_{W^{ 1, \infty} } ) \Big( |\!|V|\!|_{H^s}^2 + |\!| h |\!|_{H^{ s- 1}}^2 \Big).
\eeq
The crucial point is that this estimate only involves the $H^{s-1}$ norm of $h$. This will 
allow to conclude by using that for the first equation in \eqref{systHu}, the $H^{s- 1}$ norm 
of $h$ is controlled by the $H^s$ norm of $u$.\\

By using the commutator estimate \eqref{sob2}, we have
\begin{eqnarray*}
\mathcal{C}&  \leq & C|\!| H|\!|_{H^s} \Big(  |\!| H c(h) |\!|_{H^s } \, |\!|\nabla \cdot u |\!|_{L^\infty}
 + |\!| \nabla \big( H c(h) \big) |\!|_{L^\infty} \,  |\!|\nabla \cdot u |\!|_{H^{s-1}} \Big) \\
& \leq & C_0\big ( |\!|(h,u)|\!|_{W^{ 1, \infty}} \big) 
\Big( |\!| V |\!|_{H^s}^2 + |\!|H|\!|_{H^s}\, |\!|H c(h) |\!|_{H^s} \Big).
 \end{eqnarray*}
To estimate the last term, we use that $H= h^n$, which yields $h \p_i H = n H \p_i h$, thus
$$ \partial_{i}\Big(H c(h) \Big)= c(h) \partial_{i} H + c'(h) H \partial_{i}h 
= c(h) \partial_{i} H + \frac{1}{ n } c'(h) h \partial_{i} H.$$
Consequently, by \eqref{sob1}, \eqref{sob3}, we get
$$ |\!|H c(h) |\!|_{H^s} \leq C |\!| c(h) \nabla H |\!|_{H^{s-1}} + 
C|\!| c'(h) h \nabla  H |\!|_{H^{s-1}} 
\leq C_0 \big( |\!|(h,u)|\!|_{W^{ 1, \infty}} \big) 
\Big( |\!|H|\!|_{H^s} + |\!| h |\!|_{H^{ s- 1}} \Big),$$
and \iref{eCalpha} follows.\\

\noindent {\bf Estimate for $\boldsymbol{\BD}$.} The term $\BD$ involves derivatives 
of $u$ of order $\leq s+1$, and we shall use the energy dissipation in \iref{eHs1}. 
We prove that
\beq
\label{eDalpha}
C_1( |\!|h|\!|_{L^\ii} )\, \ep\, \mathcal{D} \leq \frac{1}{2}\, \ep\, |\!|\nabla  V |\!|^2_{H^s} + \ep \, 
C_0( |\!|h|\!|_{W^{ 1, \infty }} ) \Big( |\!|V |\!|_{H^s}^2 + |\!|\nabla h |\!|_{H^{s- 1}}^2 \Big).
\eeq

We have, on the one hand,
$$ \Big| \int_{\R^d} c'(h) \nabla h \cdot \nabla \partial^\alpha u \cdot \partial^\alpha u \Big| 
\leq C_0( |\!|h|\!|_{W^{ 1, \infty }} ) |\!|\nabla u |\!|_{H^s} \, |\! | u |\!|_{H^s}
\leq C_0( |\!|h|\!|_{W^{ 1, \infty }} ) |\!|\nabla V |\!|_{H^s} \, |\! | V |\!|_{H^s}.$$
On the other hand, for the second term (which vanishes if $n=1$), after one integration by parts 
when $|\alpha|>0$, we get
\begin{eqnarray*}
n(n-1) \Big| \int_{\R^d} \p^\alpha\big( h^{n-2} |\nabla h |^2\big) \, \p^\alpha H \Big| 
& \leq & C |\!| \nabla H|\!|_{H^s}\, |\!| h^{n-2} |\nabla h |^2 |\!|_{H^{s-1}}\\
& \leq & C_0( |\!|h|\!|_{W^{ 1, \infty }} )\, |\!| \nabla H|\!|_{H^s} |\!| \nabla h |\!|_{H^{s-1}},
\end{eqnarray*}
and if $\alpha=0$, since $H=h^n$ and $s\geq 1$,
$$ n(n-1) \Big| \int_{\R^d} h^{n-2} |\nabla h |^2 H \Big| 
= \frac{n-1}{n} \int_{\R^d} |\nabla H |^2 \leq 
C |\!| H|\!|_{H^s}^2.$$
Consequently,
$$\ep \, \mathcal{D} \leq \ep\, C_0( |\!|h|\!|_{W^{ 1, \infty }} )\, 
|\!|\nabla  V |\!|_{H^s} \Big( |\!|V |\!|_{H^s} + 
|\!|\nabla h |\!|_{H^{s- 1}} \Big)+ \ep \, C |\!|  V |\!|^2_{H^s},$$
and \iref{eDalpha} follows from the standard inequality, for $a$, $b$, $\theta>0$, 
$ab \leq \theta a^2 + \frac{b^2}{4\theta}$.\\

\noindent {\bf Estimate for $\boldsymbol{\BR}$.} We prove that
\beq
\label{eRalpha}
C_1( |\!|h|\!|_{L^\ii})\, \mathcal{R} \leq \frac{1}{2}\,\ep\, |\!| \nabla V |\!|_{H^s}^2 + 
C_0( |\!|(h,u)|\!|_{W^{ 1, \infty }} )\, |\!|V|\!|_{H^s}^2.
\eeq

By using the first equation in \eqref{hydrov} for $h$ and an integration by parts, we find, as for 
the first term in $\BD$,
\begin{eqnarray*}
\mathcal{R}^\alpha
& \leq & C_0( |\!|(h,u)|\!|_{W^{ 1, \infty }} )\, |\!|V|\!|_{H^s}^2
+ \ep\, \frac{n}{4} \int_{\R^d} c'(h) \Delta h  |\partial^\alpha u |^2 \\ 
& \leq & C_0( |\!|(h,u)|\!|_{W^{ 1, \infty }} )\, |\!|V|\!|_{H^s}^2
- \ep\, \frac{n}{4} \int_{\R^d} c'(h) \big( (\nabla h \cdot \nabla) \p^\alpha u \big) \cdot \p^\alpha u 
- \ep\, \frac{n}{4} \int_{\R^d} c''(g) |\nabla h |^2 \, |\p^\alpha u |^2 \\
& \leq & C_0( |\!|(h,u)|\!|_{W^{ 1, \infty }} )\, \Big( |\!|V|\!|_{H^s}^2 + \ep\, |\!| \nabla V|\!|_{H^s}
\, |\!|V|\!|_{H^s} \Big).
\end{eqnarray*}
Then, \iref{eDalpha} follows as above from the inequality $ab \leq \theta a^2 + \frac{b^2}{4\theta}$.\\

Summing \eqref{eCalpha}, \eqref{eDalpha} and \eqref{eRalpha}, inserting this into \eqref{eHs2} 
and simplifying by $\ep\, |\!| \nabla V |\!|_{H^s}^2$, we infer
\begin{align}
\label{eVs}
|\!| V(t)|\!|_{H^s}^2 \leq & \ C_1 \big( |\!| h(t)|\!|_{L^\ii} \big) \Big(
|\!| V(0)|\!|_{H^s}^2 \nonumber \\
& \ + \int_{0}^t C_0 ( |\!|(h,u)(\tau)|\!|_{W^{ 1, \infty}}) \Big[ |\!|V(\tau)|\!|_{H^s}^2
 + |\!| h (\tau) |\!|_{H^{s-1}}^2 + \ep |\!|\nabla h (\tau)|\!|_{H^{s-1}}^2 \Big] \, d\tau \Big).
\end{align}
To close the estimate, it remains to estimate $|\!|h |\!|_{H^{s-1}}^2$ and
$ \ep \int_{0}^t |\!| \nabla h |\!|_{H^{s-1}}^2$. We use the standard $H^{s-1}$ 
estimate for the convection diffusion equation \eqref{hydrov} which yields, as for \iref{eHs1}, 
for $|\alpha|\leq s-1$,
$$ \frac{d}{dt}\Big[ \frac{1}{2}\int_{\R^d} |\p^\alpha h|^2 \Big] + 
\ep\, \int_{\R^d} |\p^\alpha h |^2 
\leq C_0 \big( \big| \! \big| (h,u) \big| \! \big|_{W^{ 1, \ii}} \big) 
\Big( |\!|h|\!|_{H^{s-1}}^2 + |\!|h|\!|_{H^{s-1}}|\!| u |\!|_{H^{s}} \Big).$$
Summing for $|\alpha|\leq s-1$ and integrating in time, this yields
\be
\label{ehs}
\frac{1}{2}\, |\!| h(t) |\!|_{H^{s-1}}^2 + \ep\, \int_0^t |\!| \nabla h(\tau) |\!|_{H^{s-1}}^2 \, d\tau 
\leq \frac{1}{2}\,|\!|h(0)|\!|_{H^{s-1}}^2 + 
\int_0^t C_0 \big( \big| \! \big| (h,u)(\tau) \big| \! \big|_{W^{ 1, \ii}} \big) 
\Big( |\!| V(\tau) |\!|_{H^s}^2 + |\!|h(\tau)|\!|_{H^{s-1}}^2 \Big) \, d \tau.
\ee

Finally, we can combine \eqref{eVs} and \eqref{ehs}, to get 
\begin{align}
\label{euf}
|\!| V(t)|\!|_{H^s}^2 + & |\!|h(t)|\!|_{H^{s-1}}^2 \nonumber \\ \leq & 
C_0 \big( \big| \! \big| (h,u) \big| \! \big|_{L^\infty([0, t ],W^{1,\ii})} \big) 
\Big(  |\!| V(0)|\!|_{H^s}^2 + |\!|h(0)|\!|_{H^{s-1}}^2 + 
\int_{0}^t |\!| V(\tau)|\!|_{H^s}^2 + |\!|h(\tau)|\!|_{H^{s-1}}^2 \, d\tau \Big).
 \end{align}
Since $H^{s-1}$ is embedded in $W^{1,\ii}$ for $s>2+ d/2$, we easily get 
by classical continuation arguments and the Gronwall lemma that the solution of 
\eqref{hydrov} is defined on an interval of time $[0,T)$ independent of $\ep$. Finally, 
\eqref{euf} provides a uniform bound  for $(h,H,u)$ in $H^{s-1}\times H^s\times H^s$, 
which allows to prove in a classical way that $(h_\ep, u_\ep)$ converges towards a solution 
of \eqref{hydroh}. This ends the proof of the existence of solution.

To prove the uniqueness, it suffices to use the same method as above and perform an 
$L^2$ energy estimate on the system satisfied by $h_{1}- h_{2}, u_{1}-u_{2}, H_{1}- H_{2}$. 
This is left to the reader.

\subsection{WKB expansions}
\label{sectionWKB2}

\ \indent We now turn to the construction of WKB expansions up to arbitrary order. 
Let us first notice that in Theorem \ref{theoeulernh}, if the initial datum $(a_0,u_0)$ 
is in $H^\ii \times H^\ii$, then the solution $(a,u)$ is in $\BC^0([0,T],H^{s-1} \times H^s)$ 
for every $s>2+d/2$, with $T$ independent of $s>2+d/2$. In other words, the existence 
time of the maximal solution in $H^\ii \times H^\ii$ is positive. This fact follows easily 
from \iref{euf} and the Gronwall inequality (since $H^{s-1} \subset W^{1,\ii}$).

\begin{lem}
\label{WKB}
Consider $\Psi^\eps_{0}= a_{0}^\eps e^{i \varphi_0^\e / \eps}$ with 
$a_{0}^\eps \in H^\infty$, $\varphi_0^\e \in H^\infty$ and that for some 
$m \in \N$, there exists an expansion
\beq
\label{expinit} 
a_{0}^\eps = \sum_{k=0}^m \eps^k a_{0}^k + \eps^{m+1} \boldsymbol{a^\eps_{0}}, \quad 
\varphi^{\eps}_{0} = \sum_{k=0}^m \eps^k \varphi_{0}^k + \eps^{m+1} \boldsymbol{\vp_{0}^\eps}
\eeq
with $a_{0}^0 \in \mathbb{R}$, satisfying, for every $s$,
\be
\label{Hinfty}
\sup_{\eps \in (0, 1 )} \Big( \big|\!\big| \boldsymbol{a^\eps_{0}} \big|\!\big|_{H^s} 
+ \big|\!\big| \boldsymbol{\vp^\eps_{0}} \big|\!\big|_{H^s}\Big)<+ \infty.
\ee
Let us denote $T^*$ is the maximal time of 
existence of a smooth (\it{i.e.} $H^\ii \times H^\ii$) solution $(a^0,\vp^0)$ for 
\eqref{hydro} with the initial condition $(a_0^0,\vp_0^0)$. Then, there exists an 
approximate smooth solution of \eqref{NLS} on $[0,T^*)$ under the form 
$\Psi^a = a^\eps e^{i \varphi^\eps / \eps}$, with $a^\eps, \varphi^\eps \in H^\ii$ and
$a^\e$ complex-valued, solving
\be
\label{systWKB}
\left\{\begin{array}{ll}
\displaystyle{\frac{\p \vp^\e}{\p t} + f(|a^\e|^2) 
+ \frac{1}{2} |\nabla \vp^\e|^2}  = R^m_{\varphi} \\ \\
\displaystyle{\frac{\p a^\e}{\p t} 
+\big(\nabla \vp^\e \big) \cdot \nabla a^\e + \frac{a^\e}{2} \Delta \vp^\e}
- J \eps \Delta a^\eps  = R^m_{a},
\end{array}\right.
\ee
with the initial condition $\big( a^\eps, \varphi^\eps \big)_{/t=0} = \big( a_{0}^\eps, \varphi_0^\eps \big)$, 
and where, for every $s$ and $0<T<T^*$,
\be
\label{erreurreste}
\sup_{[0, T]}\Big( \big|\!\big|R^m_{a}\big|\!\big|_{H^s} + \big|\!\big| R^m_{\varphi} \big|\!\big|_{H^s} \Big) 
\leq C_{s,T} \eps^{m+2}.
\ee
Finally, for $0<T<T^*$, $a^\e$ verifies \iref{aquasireel}: $a^\e-a^0 =\BO(\e)$ in $L^\ii([0,T],W^{s,\ii})$.
\end{lem}

Note that $\Psi^a$ is indeed an approximate solution of \eqref{NLS} since
$$ i \e \frac{\p \Psi^a}{\p t} + \frac{\e^2}{2} \Delta \Psi^a 
- \Psi^a f(|\Psi^a|^2) = \Big( - i \eps R_{a}^m + a^\eps R_{\varphi}^m\Big) 
\exp \big( i \frac{\varphi^\eps}{\eps} \big).$$
By using the notation of section \ref{lineaire}, we have $R^\eps =- i\eps R_{a}^m + a^\eps R_{\vp}^m$, 
hence
\beq
\label{estReps}
\sup_{[0, T]}\big|\!\big|R^\eps \big|\!\big|_{H^s} \leq C_{s} \eps^{m+2}.
\eeq

\subsection*{Proof.}
As in \cite{G},  we look for expansions
$$ a^\eps= \sum_{k=0}^m \eps^k a^k + \eps^{M+1}a^{m+1}, \quad \quad \quad 
\varphi^\eps= \sum_{k=0}^m \eps^k \varphi^k+ \eps^{m+1}\varphi^{m+1}.$$
This yields that $(a^{0}, \varphi^0) $ solves the nonlinear system
\be
\label{nonlin0}
\left\{\begin{array}{ll}
\displaystyle{\frac{\p \vp^0}{\p t} + f(|a^0|^2) 
+ \frac{1}{2}\, |\nabla \vp^0|^2 }  = 0 \\ \\
\displaystyle{\frac{\p a^0}{\p t} 
+\big(\nabla \vp^0 \big) \cdot \nabla a^0 + \frac{a^0}{2}\, \Delta \vp^0}  =0,
\end{array}\right.
\ee
which is just \iref{aphi}, and that for $1 \leq k \leq m$, $(a^k, \varphi^k)$ solves 
the linear system
\be
\label{link}
\left\{\begin{array}{ll}
\displaystyle{\frac{\p \vp^k}{\p t} + 2 f'(|a^0|^2)(a^0, a^k) 
+ \nabla \vp^0\cdot \nabla \varphi^k }  = S^k_{\varphi } \\ \\
\displaystyle{\frac{\p a^k}{\p t} 
+\big(\nabla \vp^0 \big) \cdot \nabla a^k + \nabla a^0 \cdot \nabla \varphi^k 
+ \frac{a^0}{2}\, \Delta \vp^k + \frac{a^k}{2}\, \Delta \vp^0 }  =S^{k}_a,
\end{array}\right.
\ee
where the source terms $(S^k_{\varphi}, S^k_{a})$ depend only on $(a^j, \varphi^j)_{0\leq j \leq k-1}$, and 
$S^k_{a}$ is complex-valued.\\

We first solve \eqref{nonlin0} (that is \iref{Eulera}) with the initial condition 
$\varphi^0_{/t= 0} = \varphi^0_{0}$, $a^0_{/t=0} = a_{0}^0$. By introducing 
$u^0 \equiv \nabla \varphi^0$ and by taking the gradient of the first equation of \eqref{nonlin0}, 
we find 
\be
\label{hypernonlin}
\left\{ \begin{array}{ll}
\displaystyle{
\p_{t} a^0 +  u^0 \cdot \nabla a^0 + \frac{a^0}{2}\, \nabla \cdot u^0}  = 0 \\ \\
\displaystyle{ \p_{t} u^0 + u^0 \cdot \nabla u^0 + \nabla \Big(f \big( (a^0)^2\big) \Big)}  = 0,
\end{array} \right.
\ee
which is the compressible Euler type equation considered in the previous section. 
By using Theorem \ref{theoeulernh}, we get the existence of a smooth solution 
$(a^0, u^0) \in H^{s-1}\times H^s$ for every $s$ on $[0,T^*)$ (with $T^*$ independent 
of $s$), with $a^0$ real-valued. Finally, to get $\varphi^0$, we can
 use the same argument as in \cite{AC}.
This yields
$$\vp^0(t,x) = \vp^0_0(x) -\int_0^t \Big( f\big( (a^0)^2 \big) + \frac{1}{2}\, |u^0|^2 \Big)(\tau,x) \, d\tau.$$

We now turn to the resolution of \eqref{link}. We solve it with the initial condition 
$\big( \varphi^k, a^k\big)_{/t=0} = \big( \varphi_0^k, a_0^k \big)$. By introducing again 
$u^k \equiv \nabla \varphi^k$, we can take the gradient in the first line of \eqref{link} to get
\be
\label{hyperak}
\left\{ \begin{array}{ll}
\displaystyle{  \p_{t} a^k + u^0 \cdot \nabla a^k +  \frac{a^0}{2} \nabla \cdot u^k + 
u^k \cdot \nabla a^0  + \frac{a^k}{2}  \nabla \cdot u^0}  = S^k_{a}, \\ \\
\displaystyle{  \p_{t } u^k + u^0 \cdot \nabla u^k  + \nabla\big( a^0, f'( (a^0)^2) a^k \big) 
+ u^k \cdot \nabla u^0 } 
= \nabla S^k_{\varphi}.
\end{array}
\right.
\ee
Again, since $f'\big( (a^0)^2\big)$ can vanish, the symmetrization of this linear hyperbolic system 
requires some care. We thus  set
$$F^k (t,x) \equiv 
\left\{ \begin{array}{ll}\sqrt{2} \,\big( f'((a^0)^2)\big)^{\frac{1}{2}} a^k \quad  \quad \,
 \mbox{ if $n$ is odd} \\
 \sqrt{2} \, a^{0}\big( {f'((a^0)^2) \over (a^{0})^2}\big)^{\frac{1}{2}} a^k \quad \quad
 \mbox{ if $n$ is even}  \end{array} \right.$$
 Note that in both cases, we have
 $$  F^k(t,x) = \sqrt{2} \, g(a^0) a^k$$
 with $g$ smooth.  Indeed,  by using that we can write
  $f'(\rho)= \rho^{n-1} \hat{f}(\rho)$ with $f$ smooth and positive, we have
   in both cases :
 \beq
 \label{formeg}
 g(a^0)=    (a^0)^{n-1}\,\big( \hat{f}((a^0)^2)\big)^{\frac{1}{2}}.
\eeq
This is the  natural generalization of the change of unknown used in \cite{AC}.
 Then, thanks 
to the equation on $a^0$, we get for $(F^k, u^k)$ the system
$$\left\{ \begin{array}{ll}
\displaystyle{  \p_{t} F^k + u^0 \cdot \nabla F^k + \frac{1}{\sqrt{2}}\, a^0 g(a^0) \nabla \cdot u^k
\ \ + \sqrt{2}\, g(a^0) \, u^k \cdot \nabla a^0 + \frac{F^k}{2} \Big( 1 + \frac{a^0 g'(a^0)}{ g(a^0)} \Big) 
\nabla \cdot u^0 }
 = \sqrt{2} g(a^0) S^k_{a}, \\ \\
\displaystyle{ \p_{t } u^k + u^0 \cdot \nabla u^k \ + \frac{1}{\sqrt{2}}\, \nabla\big (a^0 g(a^0), F^k \big) 
+  u^k \cdot \nabla u^0 }  = \nabla S^k_{\varphi}.
\end{array}
\right.$$
Note that the coefficient $ \frac{a^0 g'(a^0)}{ g(a^0) }$ is smooth even 
when $a^0$ vanishes since $g$ is under the form \eqref{formeg}. We have obtained 
a linear symmetric hyperbolic system with a zero order term and a source term $\BS^k$ 
depending only on $(a^j,\vp^j)$ for $0\leq j <k$ under the form
$$ \partial_{t} U^k + \sum_{j= 1}^d A^j (t,x) \partial_{j}U^k + L(t,x) U^k= \BS^k, \quad 
U^k = \left( \begin{array}{c} F^k\\ u^k \end{array} \right),$$
where $ A^j(t,x)$ are smooth, real and symmetric and the matrix $L$ is smooth. By the classical 
theory, there exists, on $[0,T^*)$, a smooth solution $(F^k, u^k)$ in $H^\ii\times H^\ii$ of this system. 
 Once $u^k$ is built, we get $a^k$ by solving the transport equation for $a^k$ 
which is given by the first line of \eqref{hyperak}. Finally, we deduce the phase $\vp^k$ by integrating 
in time the first line of \eqref{link}. We obtain
$$ \varphi^k (t,x) = -\int_{0}^t \Bigl( 2 f'(|a^0|^2)(a^{0}, a^k) + \nabla \vp^0\cdot u^k 
- S^k_{\varphi } \Bigr)(\tau,x) d\tau.$$ 

Finally, we choose in a similar way $(a^{m+1}, \varphi^{m+1})$ that solve \eqref{link} 
with the initial condition $\big( a^{m+1},\vp^{m+1}\big)_{/t=0} = 
\big( \boldsymbol{a^\eps_0}, \boldsymbol{\varphi^{\e}_0}\big)$. Because 
of the assumption \eqref{Hinfty}, we find that they are also uniformly bounded in $H^{s-1}\times H^s$ 
with respect to $\eps$. This concludes the proof of Lemma \ref{WKB}.\b

\section{Nonlinear stability}
\label{nonlineaire}

\ \indent In this section, we give the proof of Theorem \ref{Conver}. We shall actually prove 
directly a more precise version which states the existence of a WKB expansion to any order.

\begin{theo}
\label{WKBstab}
Consider $\Psi^\eps_{0}= a_{0}^\eps e^{i \varphi_0^\eps / \eps}$ with $a_{0}^\eps \in H^\infty$, 
$\varphi_0^\eps \in H^\infty$ and that for some $m \in \N$, there exists an expansion 
\eqref{expinit} as in Lemma \ref{WKB}. We assume $(\mathcal{A})$ and let $(a^\eps, \varphi^\eps)$ 
be the smooth approximate solution given by Lemma \ref{WKB} which is smooth on $[0, T^*)$. Then,\\

\noindent $\bullet$ if $m=0$,  there exists $\eps_{0}>0$ and $T\in (0,T^*)$ such that for every 
$ \eps \in (0, \eps_{0}]$, the solution of \eqref{NLS} with initial data $\Psi_{0}^\eps$ remains 
smooth on $[0, T]$ and satisfies for every $s \in \mathbb{N}$, the estimate
$$ \bigg|\! \bigg| \Psi^\e \exp\big( -\frac{i}{\e} \vp^\e \big) - a^\e \bigg|\! \bigg|_{L^{\ii}([0,T],H^s)} 
\leq C_{s} \eps.$$

\noindent $\bullet$ if $m\geq 1$, for every $T\in (0,T^*)$, there exists $\eps_{0}(T)>0$ such 
that for every $\eps \in (0, \eps_{0}(T)]$, the solution of \eqref{NLS} with initial data $\Psi_{0}^\eps$ 
remains smooth on $[0, T]$ and satisfies for every $s \in \mathbb{N}$, the estimate
$$ \bigg|\! \bigg| \Psi^\e \exp\big( -\frac{i}{\e} \vp^\e \big) - a^\e 
\bigg|\! \bigg|_{L^{\ii}([0,T],H^s)} \leq C_{s,T} \eps^{m+1}.$$ 
\end{theo}

Note that Theorem \ref{Conver} is actually the special case $m=0$ in Theorem \ref{WKBstab}.

\subsection*{Proof of Theorem \ref{WKBstab}}

\ \indent Let $s > d/2$. We take $(a^\eps, \varphi^\eps)$ the approximate solutions given by 
Lemma \ref{WKB} and look for  the solution of \eqref{NLS} under the 
form $\Psi^\eps =( a^\eps + w) e^{i \varphi^\eps / \eps} $. We get for $w$ the 
equation \eqref{NLSwlin} with $F^\eps$ given by \eqref{Feps} and the initial condition 
$w_{/t=0}=0$. For $s>d/2$, and every $\eps>0$, this semilinear  equation is locally well-posed in 
$H^s$: we get very easily that there exists for some $T^\eps >0$ a unique maximal 
solution $w \in \mathcal{C}([0, T^\eps ), H^s) $ of \eqref{NLSwlin} (see \cite{Cazenave} 
for example). We shall prove that $T^\eps$ is bounded from below by some  $T>0$ if $m=0$, and  that 
$T^\e \geq T$ for every $T \in (0, T^*)$ for $\eps$ sufficiently small if $m\geq 1$. Let us define
$$ \tau^\e \equiv \sup \big\{ \tau \in (0,T^\e), \ \forall t\in [0,\tau], \ 
2 N^\eps_s\big( w(t) \big) \leq \e^{2m+4} \big\}.$$
Note that $\tau^\eps>0$ since $w(0)=0$ and that by Sobolev embedding, we have, for 
$t \leq \tau^\eps$,
$$ \big|\! \big| w(t) \big|\! \big|_{L^\infty}^2 \leq 
K^2 \eps^{-2} N_{s}^\eps(w(t))\leq K^2 \eps^{2m+2} \leq K^2,$$
for some $K$ independent of $\eps$.\\

We will apply Theorem \ref{estiHS} with $F^\eps$ given by \eqref{Feps}. To estimate 
$F^\eps$, we use the following lemma:

\begin{lem}
\label{estirhs}
Let $R>0$, $s>d/2$ and $w$ such that $\big|\!\big| w \big|\!\big|_{L^\ii} \leq R$, and $F^\e$ given 
by \eqref{Feps}. Then, for a constant $C$ depending only on 
$\big|\!\big| a^\e(t) \big|\!\big|_{W^{s+2,\ii}}$ and $R$, we have
$$ \big|\!\big| F^\e \big|\!\big|_{H^s}^2 
+ \frac{1}{\e^2} \big|\!\big| {\rm Im} F^\e \big|\!\big|_{H^{s-1}}^2 
\leq C \e^{2m+4} + C\eps^{2m} N^\eps_{s}(w) + C \bigg[ \frac{N^\eps_s(w)}{\e^4} + 
\bigg( \frac{N^\eps_s(w)}{\e^4} \bigg)^2 \bigg] N^\eps_s(w).$$
\end{lem}

We postpone the proof of Lemma \ref{estirhs} to the end of the section. We can first 
easily end the proof of Theorem \ref{WKBstab}. Notice first that, by definition of $\Psi^a$, 
we have
$$ R_a = R_a^m + i\e \Delta a^\e = \BO_{H^k}(\e^{m+1}) + \BO_{H^k}(\e) = \BO_{H^k}(\e),$$
for every $k$, uniformly for $0\leq t \leq T$, hence
$$ \frac{ 1 }{\e} \big|\!\big| R_a (t) \big|\!\big|_{W^{s-1,\ii}} \leq C.$$
Applying Theorem \ref{estiHS} and Lemma \ref{estirhs} with 
$R\equiv K$, we infer that for $0\leq t \leq \tau^\e$,
$$ \frac{d}{dt} N^\eps_s \big( w(t) \big) \leq C \e^{2m+4} + C  \eps^{2m } N^\eps_s \big( w(t) \big),$$
which gives immediately, since $w_{/t=0}=0$, that
$$N^\eps_s \big( w(t) \big) \leq C \e^{2m+4} \, \Big( e^{C \eps^{2m}t } - 1\Big) \leq \frac{1}{2}\, \e^{2m +4}$$
in the following cases:
\begin{itemize}
\item for $m=0$, $0\leq t\leq T $ with $0<T< T^*$ sufficiently small independent of $\eps$,
\item for $m \geq 1$, $T\in (0,T^*)$ is arbitrary, $0\leq t \leq T$ and $ \eps \leq \eps_{0}(T)$ 
with $\eps_{0}(T)$ sufficiently small.
\end{itemize}
As a consequence,
$ \tau^\e \geq T$
as desired and
$$ \big|\!\big| w \big|\!\big|_{L^\ii ([0,T], H^s(\R^d))} \leq C_{s,T} \e^{m+1}.$$

It remains to prove Lemma \ref{estirhs}.

\subsubsection*{Proof of Lemma \ref{estirhs}}

\ \indent
We recall that $F^\eps$ is given by
$$ F^\eps = R^\eps + Q(w) = R^\e + (a^\eps + w ) \Big( f( |a^\eps + w|^2) - f(|a^\eps |^2) \Big) - 
2 (w,a^\eps) f'(|a^\eps|^2) a^\eps.$$
As a first try, we could use the rough estimate 
$$ Q^\eps(w) = \BO( |w|^2)\quad \quad \quad {\rm as} \ \ w \to 0,$$
which would lead to
$$\big|\!\big|  Q^\e \big|\!\big|_{H^s}^2 
+ \frac{1}{\e^2} \big|\!\big| {\rm Im}\, Q^\e \big|\!\big|_{H^{s-1}}^2 \leq 
\frac{C}{\e^2} \big|\!\big| w \big|\!\big|_{H^s}^4 \leq \frac{C}{\e^6} N^\eps_s(w)^2,$$
which does not allow to conclude in the proof of Theorem \ref{WKBstab} for $m=0$ 
and does not give a sharp result for the existence time if $m=1$. To get the refined 
estimate of Lemma \ref{estirhs}, the idea is then to use a Taylor expansion 
for $Q^\e$ w.r.t. $w$ up to second order, and write
$$Q^\eps(w) = |w|^2 f'(|a^\e|^2) a^\e + 2 f'(|a^\e|^2) (w,a^\e) w  + 
2 a^\eps f''(|a^\eps |^2) (w,a^\eps)^2 + G^\e(x,w),$$
so that for fixed $x$, we have as $w \to 0$,
$$ G^\e(x,w) = \BO\big(|w|^3 \big).$$
We turn now to estimate each term in $F^\e$.\\

\noindent {\it Estimate for} $\displaystyle {R^\eps = i \eps R_{a}^m - R_{ \varphi}^m a^\eps}$. 
Thanks to \eqref{estReps}, we have
$$\big|\!\big| R^\e \big|\!\big|_{H^s}^2 \leq C \e^{2m+ 4}.$$
Moreover, since $R_{\varphi}^m$ is real-valued and since, from \iref{aquasireel}, 
${\rm Im}\, a^\e = \BO_{W^{s,\ii}}(\e)$, we also have
$$\frac{1}{\e^2} \big|\!\big| {\rm Im}\, R^\e \big|\!\big|_{H^{s-1}}^2 \leq C\e^{2m + 4}$$
thanks to \eqref{erreurreste}. We have thus proven that
\be
\nonumber
\big|\!\big| R^\eps \big|\!\big|_{H^s}^2 
+ \frac{1}{\e^2} \big|\!\big| {\rm Im}\, R^\e \big|\!\big|_{H^{s-1}}^2 \leq C \e^{2m+ 4 }.
\ee
\bigskip

\noindent {\it Estimate for} $\displaystyle{G^\e(x,w)}$. The estimate relies on 
Lemma \ref{composition} in the appendix. Indeed, it is clear from the Taylor formula that $G^\e$ 
may be written under the form
$$ \big( {\rm Re}\, w \big)^2 h_{11}\big( x,w(x) \big) + 
\big( {\rm Re}\, w \big) \big( {\rm Im}\, w \big) h_{12}\big( x,w(x) \big) +
\big( {\rm Im}\, w \big)^2 h_{22}\big( x,w(x) \big),$$
where $h_{11}$, $h_{12}$, $h_{22} : \R^d \times \C \to \C$ are of class $\BC^\ii$ and 
$\forall x\in \R^d$, $h_{11}(x,0) = h_{12}(x,0) = h_{22}(x,0) = 0$. Moreover, $h_{11}$, 
$h_{12}$ and $h_{22}$ verify the hypothesis of Lemma \ref{composition} in the Appendix since  $a^\eps \in L^\infty([0, T], H^s).$ 
As a consequence, if $\big|\!\big| w \big|\!\big|_{L^\ii} \leq R$,
$$\big|\!\big|  G^\e \big|\!\big|_{H^s}^2 \leq C \big|\!\big| w \big|\!\big|_{H^s}^3,$$
which implies
\be
\nonumber
\big|\!\big| G^\e \big( x,w(x) \big) \big|\!\big|_{H^s}^2 + 
\frac{1}{\e^2} \big|\!\big| {\rm Im}\, G^\e\big( x,w(x) \big) \big|\!\big|_{H^{s-1}}^2 \leq 
\frac{2}{\e^2} \big|\!\big| G^\e\big( x,w(x) \big) \big|\!\big|_{H^s}^2 \leq \frac{C}{\e^8} N^\eps_s(w)^3.
\ee

\bigskip

The estimate for the quadratic terms in $Q^\e(w)$ will rely crucially on the fact 
that $a^\e$ is {\it real} to first order and  that $(w,a^\e)$ is estimated in 
$H^{s-1}$ by $N^\e_s(w)$ and not just by $\e^{-2} N^\e_s(w)$.\\

\noindent {\it Estimate for} $\displaystyle{F^\e_1 \equiv |w|^2 f'(|a^\e|^2) a^\e}$. We have
$$ \big|\!\big| F^\e_1 \big|\!\big|_{H^s}^2 \leq \frac{C}{\e^4} N^\eps_s(w)^2,$$
and in view of \iref{aquasireel}, Im $a^\e = \BO_{W^{s,\ii}}(\e)$, thus
$$ \frac{1}{\e^2} \big|\!\big| {\rm Im}\, F^\e_1 \big|\!\big|_{H^{s-1}}^2 
\leq C \big|\!\big| |w|^2 \big|\!\big|_{H^{s-1}}^2 \leq \frac{C}{\e^4} N^\eps_s(w)^2.$$

\bigskip

\noindent {\it Estimate for} $\displaystyle{F^\e_2 \equiv 2 f'(|a^\e|^2) (w,a^\e) w}$. We 
begin with the rough estimate
$$ \big|\!\big| F^\e_2 \big|\!\big|_{H^s}^2 \leq \frac{C}{\e^4} N^\eps_s(w)^2.$$
Moreover, one has
\be
\label{prod}
\big|\!\big| f'(|a^\e|^2) (w,a^\e) \big|\!\big|_{H^{s-1}}^2 \leq C N^\e_s(w).
\ee
Indeed, let $\mu \in \N^d$ with $|\mu| \leq s-1$. Then,
$$ \p^\mu \big( f'(|a^\e|^2) (w,a^\e) \big) = \sum_{ \alpha +\beta + \lambda =\mu}\ * \ 
\p^\lambda \big[ f'(|a^\e|^2) \big] \big(\p^\alpha w, \p^\beta a^\e\big),$$
where $*$ is a coefficient depending only on $\alpha$, $\beta$ and $\lambda$. Since 
$|\mu|\leq s-1$, the terms $(\p^\alpha w, \p^\beta a^\e)$ are bounded in $L^2$ by 
$\Sigma(w)^{\frac{1}{2}} + \eps |\!|w|\!|_{H^{s-2}}$ as soon as $|\alpha |\leq s-2$.
 The  term in the sum 
with $|\alpha | = s-1$ (hence $\mu=\alpha$ and $\beta =\lambda =0$) is 
$f'(|a^\e|^2)\big(\p^\mu w, a^\e\big)$ is bounded in $L^2$ by $N^\e(\p^\mu w)$. Hence, 
\iref{prod} follows.

As a consequence, by \iref{sob1} and Sobolev embedding, we obtain
$$\big|\!\big| f'(|a^\e|^2) (w,a^\e) w \big|\!\big|_{H^{s-1}} \leq 
C_s \big|\!\big| w \big|\!\big|_{L^\ii} \Big( \big|\!\big| f'(|a^\e|^2) (w,a^\e) \big|\!\big|_{H^{s-1}} + \big|\!\big| w \big|\!\big|_{H^{s-1}} \Big) \leq \frac{C}{\e^2} N^\e_s(w).$$
Consequently,
$$\big|\!\big| F^\e_2 \big|\!\big|_{H^s}^2 + 
\frac{1}{\e^2} \big|\!\big| {\rm Im}\, F^\e_2 \big|\!\big|_{H^{s-1}}^2 
\leq \frac{C}{\eps^4} N^\eps_{s}(w)^2.$$

\bigskip

\noindent {\it Estimate for} $\displaystyle{F^\e_3 \equiv 2 a^\eps  f''(|a^\e|^2) (w,a^\e)^2}$. 
We find as for $F^\e_1$
$$ \big|\!\big| F^\e_3 \big|\!\big|_{H^s}^2 \leq \frac{C}{\e^4} N^\eps_s(w)^2,$$
and once again in view of \iref{aquasireel},
$$\frac{1}{\e^2} \big|\!\big| {\rm Im}\, F^\e_3 \big|\!\big|_{H^{s-1}}^2 
\leq C \big|\!\big| w \big|\!\big|_{H^{s-1}}^4 \leq \frac{C}{\e^4} N^\eps_s(w)^2.$$

We conclude the proof of Lemma \ref{estirhs} summing these estimates. \b

\section{Geometric optics in a half-space}
\label{half}

\ \indent In this section, we consider the Gross-Pitaevskii equation in a half-space 
in dimension $d\leq 3$
\beq
\label{GP+}
GP(\Psi^\eps) \equiv i \eps \partial_{t}\Psi^\eps + \frac{\eps^2}{2} \Delta \Psi^\eps - 
\Psi^\eps( |\Psi^\eps |^2 - 1 ) = 0, \quad  x \in  \R^d_+ \equiv \R^{d-1} \times (0,+\ii).
\eeq
We consider the Neumann 
boundary condition \eqref{dir+} on the boundary and the condition \eqref{inf+} at infinity, 
that is
$$ \frac{\p \Psi^\e}{\p n}_{/\p \R^d_+} = \frac{\p \Psi^\e}{\p z}_{/ z=0} =0 \quad \quad 
{\rm and} \quad \quad \exp\Big( - \frac{i}{2\e} |u^\ii|^2 \,t - \frac{i}{\e} u^\ii \cdot x \Big) \Psi^\e 
\to 1 \quad |x| \to +\ii$$
by using the notation
 $\ x =(y,z) \in \R^{d-1} \times (0,+\ii)$.

\subsection{Construction of the WKB expansion}

\ \indent In this section, we shall consider a smooth solution $(a,u)$, with $a$ real-valued, of
\be
\label{eul+}
\left\{ \begin{array}{ll}
\displaystyle{ \p_t a + \, u \cdot \nabla a + \frac{1}{2}\, a  \, \nabla \cdot u}  = 0 \\ \\
\displaystyle{ \p_t u + u \cdot \nabla u + \nabla( a^2) }  = 0,
\end{array}
\right.
\ee
with the boundary condition $u_{d}(t,y, 0)= 0$ and the condition at infinity
$$ u(t,x) \rightarrow u^\infty, \quad \quad a(t,x) \rightarrow 1 \quad \quad {\rm when} \ \ 
|x | \rightarrow + \infty.$$
Since we look for $a$ real-valued, the resolution of this system is made in \cite{LZ} (Theorem 2). 
Given $s\in \N^*$, if the initial datum $a_0$ is positive and $(a_0 -1,u_0-u^\ii)\in H^s$, and 
under some compatibility conditions for $(a_0,u_0)$ on the boundary $\p \R^d_+$ of sufficiently high 
order on the initial data, there exists $T_0\in (0,+\ii)$ and a solution $(a,u)$ on $[0, T_0]$ with 
$(a-1,u-u^\ii)\in \BC^0([0,T_0],H^s) \cap \BC^1([0,T_0],H^{s-1})$, such that
\beq
\label{apos}
a(t,x) \geq \alpha >0, \quad \forall t \in [0, T_0], \, \forall x \in \overline{\R^d_+}.
\eeq
for some $\alpha >0$. We also define the phase $\varphi$ by
$$\vp(t,x) \equiv \vp_0(x) - \int_{0}^t \Big( \frac{1}{2} | \nabla \vp |^2 + |a|^2 - 1\Big) (\tau,x)\ d\tau.$$
In view of the condition \iref{inf+} at infinity, $\vp$ is not in $H^s$ but 
$\vp(t, .) - u^\ii \cdot x - \frac{t}{2} |u^\ii|^2 \in H^s$.\\

The aim of this subsection, is to prove the existence of WKB expansion (which involves 
boundary layers since the solution of \eqref{eul+} does not match the Neumann boundary 
condition \eqref{dir+}) up to arbitrary orders for \eqref{GP+}, \eqref{dir+}, \eqref{inf+} 
starting from a smooth $(a,u)$ which verifies \eqref{apos}.

We define the set of  boundary layer profiles $\mathcal{S}_{exp}$ as
$$ \mathcal{S}_{exp} = \Big\{ A(t,y,Z) \in H^\infty(\R_{+ }\times \R^{d-1} \times \R_{+}), 
\quad \forall k, \, \alpha, \, l, \quad \exists \gamma>0, \quad | \p_t^k \p_y^\alpha \p_{Z}^l A| \leq 
C_{k, \alpha, l} \exp(-\gamma Z) \Big\}.$$

\begin{lem}
\label{WKB+}
Let $s\in \N$ be fixed. For every $m\in \N^*$, there exists a smooth approximate solution 
$\Psi^{a,m}= a^\eps e^{i \frac{\varphi^\eps}{\eps }}$ on $[0, T_m]$ of \eqref{GP+}, with 
the Neumann condition \iref{dir+} and the condition \iref{inf+} at infinity, such that
\beq
\label{erreur+}
GP(\Psi^{a, m})= \eps^{m} R^\eps e^{i \frac{\varphi^\eps}{\eps}},
\eeq
where $R^\eps$ can be written under the form
\beq
\label{formR+}
R^\eps = - a^\e \Big( R^{int,m}_{\varphi}(t,x)+ R^{\flat,m}_{\varphi}(t,y,\frac{z}{\e}) \Big) + 
i \Big( \eps R_{a}^{int, m}(t, x) + R_{a}^{\flat,m}(t, y, \frac{z}{\e}) \big),
\eeq
with $R^{int, m}_\vp$, $R^{int,m}_a$ smooth and uniformly bounded in $H^s$ and 
$R_{a}^{\flat,m}(t,y,Z)$, $R_{\varphi}^{\flat,m}(t,y,Z) \in \mathcal{S}_{exp}$. Moreover, 
$\Psi^{a,m}$ verifies \eqref{dir+}, \eqref{inf+}, $a^\eps$ is real-valued and $a^\eps$, $\varphi^\eps$ 
have smooth expansions under the form
\begin{eqnarray}
\label{exp+}
a^\eps & = & a + \sum_{k= 1}^{m-1}\eps^k \Bigl( a^k(t,x) + A^k (t,y, \frac{z}{\eps}) \Big) + 
\eps^m A^m(t,y,\frac{z}{\eps}), \\
\label{exp+2}
\varphi^\eps & =& \varphi + \sum_{k=1}^{m-1} \eps^k\Big( \varphi^k(t,x) + \Phi^k(t,y, \frac{z}{\eps} ) \Big)
+ \eps^m \Phi^m(t, y, \frac{z}{\eps}).
\end{eqnarray} 
The boundary layer profiles $A^k (t,y,Z)$, $\Phi^k(t, y, Z)$ belong to $\mathcal{S}_{exp}$ 
and are such that
$$ \partial_{Z}A^1 (t,y,0) = - \p_{z} a(t,y,0), \quad \quad \quad \quad 
\partial_{Z}\Phi^1 (t,y,0) = - \p_{z} \vp(t,y,0),$$
\beq
\label{bordk}
\partial_{Z}A^k (t,y,0) = - \partial_{z} a^{k- 1}(t,y,0), \quad  \quad 
\partial_{Z}\Phi^k (t,y,0) = - \partial_{z} \varphi^{k- 1}(t,y,0)\quad \quad \quad \forall 2 \leq k \leq m.
\eeq
\end{lem}

\subsection*{Proof :} Since $\Psi^{a,m} = a^\eps \exp\big( i \frac{\varphi^\eps}{\eps} \big)$, 
we want to solve approximately
\beq
\label{phia+}
- a^\eps\Big( \partial_{t} \varphi^\eps + \frac{1}{2} | \nabla \varphi^\eps |^2 + | a^\eps |^2 - 1  \Big)
+ i \eps \Big( \partial_{t} a^\eps + \nabla \varphi^\eps \cdot \nabla a^\eps + 
\frac{1}{2}\, a^\eps \Delta \varphi^\eps \Big) + \frac{\eps^2}{2} \Delta a^\eps = 0.
\eeq
Since, in this section, we are looking for $a^\eps$ real-valued, we can split the system \eqref{phia+} 
into
\be
\label{eulkort}
\left\{ \begin{array}{ll}
\displaystyle{\p_t a^\eps + \nabla \vp^\eps \cdot \nabla a^\eps + \frac{1}{2}\, a^\e \Delta \vp^\eps}  = 0 \\ \\
\displaystyle{ \p_t \vp^\eps + \frac{1}{2} | \nabla \vp^\eps |^2 + (a^\e)^2 - 1}  = 
\displaystyle{\frac{\eps^2}{2}\, \frac{\Delta a^\eps }{a^\eps}}
\end{array}\right. \quad \quad \quad {\rm for} \ \ t \geq 0, \ \ x \in \R^d_+.
\ee
Note that in this section, the division by $a^\eps$ in the right-hand side of the second equation 
of \eqref{eulkort} is not a problem since $a^0=a$ verifies \eqref{apos} and hence does not vanish.

We thus plug the expansions \eqref{exp+}, \eqref{exp+2} in \eqref{eulkort} and 
we cancel the powers of $\eps$. To separate interior and boundary layer terms, we use 
the general theory of \cite{Grenier-Gues}. In particular, we use that for every 
function $f$, we have the expansion
$$f\Big( u(t,x)+ V(t,y, z/\eps) \Big)= f\Big(u(t,x)\Big) + f\Big(u(t,y,0) + V(t,y,z/\eps)\Big)- 
f\Big(u(t,y,0)\Big) + \eps \BR,$$
where $\BR\in \BS_{exp}$. This yields that the boundary layer part of $f\big(u(t,x)+ V(t,y,z/\eps) \big)$ 
is given by $f\big(u(t,y,0) + V(t,y,z/\eps)\big)- f\big(u(t,y,0)\big)$. In the following, we use the notation 
$W_{b}= W(t,y,0)$ for every $W(t,x)$. At first, the $\eps^{-1} $ term in the equation only gives 
$$ a_{b} \partial_{ZZ} \Phi^1 = 0 $$
and hence we have $\Phi^1 = 0$, since $a_b \geq \alpha >0$ and $\Phi^1 \in \BS_{exp}$. Note that 
this is coherent with the fact that $u_{d}(t,y,0) = \partial_{z} \varphi_{b}=0$ so that we do not 
need a boundary layer to correct the boundary condition. The $\eps^0$ term gives
\be
\label{0+}
\left\{ \begin{array}{ll}
\displaystyle{\partial_{t} \varphi + \frac{1}{2} |\nabla \varphi |^2 + a^2-1}  = 0 \\ \\
\displaystyle{\partial_{t} a + \nabla \vp \cdot \nabla a  + \frac{1}{2}\, a\, \Delta \varphi }  =0
\end{array}\right. \quad \quad {\rm for} \ \ t\geq 0, \ \ \ x \in \R^d_+
\ee
for the interior part (which is the expected equation) and for the boundary layer terms, for 
$(t,y)\in \R^+ \times \R^{d-1}$,
\beq
\label{2b+}
a_{b} \partial_{ZZ} \Phi^2 = - (\p_{z} \varphi)_b\, \partial_{Z} A^1 = 0 
\quad \quad {\rm for} \ \ Z>0,
\eeq
since $(\p_{z}\vp)_{b}= u_{d}(t,y,0)=0$. Consequently, we also find $\Phi^2 = 0$. Next, the 
order $\eps$ gives
$$\left\{ \begin{array}{ll}
\displaystyle{\p_{t}a^1 + \nabla \varphi \cdot \nabla a^1 + \nabla \varphi^1 \cdot \nabla a + 
\frac{1}{2}( a \Delta \varphi^1 + a^1 \Delta \varphi)}  = 0 \\ \\
\displaystyle{\partial_{t} \varphi^1 +  2 a\, a^1 + \nabla \varphi \cdot \nabla \varphi^1 }  =0
\end{array}\right. \quad \quad {\rm for} \ \ t \geq 0, \ \ \ x \in \R^d_+$$
in the interior and for the boundary layer terms
\be
\label{1b}
\left\{ \begin{array}{ll}
\displaystyle{\frac{1}{2} \partial_{ZZ} A^1}  = 
A^1 \Big( \partial_{t} \varphi + \frac{1}{2} | \nabla \varphi|^2 + a^2 -1 \Big)_b + 
2 a_{b}^2 A^1 = 2 a_{b}^2 A^1\\ \\
\displaystyle{a_b \partial_{ZZ} \Phi^3}  = G^1
\end{array}\right. \quad \quad \quad {\rm for} \ \ Z>0,
\ee
where $G^1 \in \mathcal{S}_{exp}$ depends only on $(a,A^1, a^1)$ and 
$(\varphi, \vp^1)$. Consequently, the boundary layer $A^1$ is given by
$$ A^1 \equiv \frac{ (\p_{z} a )_{b}}{2 a_b} e^{- 2a_{b} Z}$$
in order to match \eqref{bordk}. Finally, the $\eps^k$, $k \geq 2$ term gives
\be
\label{eqlatea+phi+}
\left\{ \begin{array}{ll}
\displaystyle{\p_t \vp^k + 2 a \, a^k + \nabla \varphi \cdot \nabla \varphi^k}  = S^k_{\varphi}\\ \\
\displaystyle{\partial_{t} a^{k} + \nabla \varphi \cdot \nabla a^{k} + \nabla a \cdot \nabla \varphi^k  + 
\frac{a}{2} \Delta \varphi^k  + \frac{a^k}{2} \Delta \varphi}  = S^k_{a}
\end{array}\right. \quad \quad \quad {\rm for} \ \ t \geq 0, \ \ x \in \R^d_+,
\ee
and
\be
\label{eqlateA+Phi+}
\left\{ \begin{array}{ll}
\displaystyle{\partial_{ZZ} A^{k}}  = 4 a_{b}^2 A^k+ F^k \\ \\
\displaystyle{\partial_{ZZ} \Phi^{k} }  = G^k
\end{array}\right. \quad \quad \quad {\rm for} \ \ Z > 0,
\ee
where $S^k_{\varphi}$ and $S^k_{a}$ depend only on $(a,\vp)$, $(a^j,\varphi^j)_{1\leq j \leq k- 1}$, 
$F^k \in \mathcal{S}_{exp}$ depend only on $(a,\vp)$, $(a^j,\vp^j,A^j,\Phi^j)_{1\leq j \leq k-1}$ 
and $(a^k,\vp^k,\Phi^k)$, and $G^k \in \mathcal{S}_{exp}$ depend on $(a,\vp)$, 
$(a^j, \vp^j,A^j,\Phi^j)_{1\leq j \leq k-1}$. Therefore, if we want to solve by induction these 
equations, one has to determine first $\Phi^k$, then $(a^k,\vp^k)$ and finally $A^k$.

To solve the cascade of equations by induction, we first find $(a^1, \varphi^1)$. As before, we 
notice that $(a^1, u^1 \equiv \nabla \varphi^1)$ solves a symmetrizable hyperbolic system (there 
is no problem with the vacuum since we are in the same situation as in \cite{G}). Since the 
condition at infinity is already absorbed by $(a, \varphi)$, one can look for $(a^1, u^1)$ in $H^s$. 
Moreover, we solve the system in $\R^d_+$ with the boundary condition $ u^1_{d}(t,y,0)=0$ which is 
needed in order to match \eqref{bordk} since we have already found that $\Phi^2 = 0$. The existence 
of a smooth solution for this linear system with the boundary condition $u^1_{d}(t, y, 0) = 0$ which 
is maximal dissipative can be obtained by the classical theory \cite{Rauch}. Then, one finds 
$\vp^1$ by the formula
$$ \vp^1(t,x) = - \int_0^t \big( 2a\, a^1 + u \cdot u^1\big) (\tau,x)\ d\tau.$$
Furthermore, since $F^2 \in \BS_{exp}$ and $a_b \geq \alpha >0$, the first equation in 
\iref{eqlateA+Phi+} (with $k=2$) has a unique solution $A^2 \in \BS_{exp}$. We have therefore 
found $(a^1,A^1,\vp^1,\Phi^1,A^2,\Phi^2)$.

We now proceed by induction. Assume that, for some $m\geq 2$, we have determined 
$(a^{j}, \vp^{j})_{1\leq j \leq m-1}$ and $(A^{j}, \Phi^{j})_{1\leq j \leq m}$. Then, we 
wish to solve \iref{eqlatea+phi+} and \iref{eqlateA+Phi+} with $k=m+1$. Since $G^{m+1}$ is already 
determined and $G^{m+1}\in \BS_{exp}$, the differential equation $\partial_{ZZ} \Phi^{m+1}= G^{m+1}$ has 
a unique solution in $\BS_{exp}$ and
$$ \p_Z \Phi^{m+1}(t, y, Z) = - \int_{Z}^{+ \infty} \frac{G^{m+1}(t,y,\zeta)}{a_{b}(t,y)} \, d\zeta.$$
This determines the boundary condition for $u^{m+1} \equiv \nabla \varphi^{m+1}$. Indeed, to 
match \eqref{bordk} we shall need to impose
\be
\label{transcondbord}
u^{m+1}_{d}(t,y,0) = (\p_z \vp^{m+1})(t,y,0) = - (\p_Z \Phi^{m+1})(t,y,0) = 
\int_{0}^{+ \infty} \frac{G^{m-1}(t,y,\zeta)}{a_{b}(t,y)} \, d\zeta,
\ee
which is non-zero in general. We then solve \iref{eqlateA+Phi+} in the following way: 
$(a^{m+1}, u^{m+1} \equiv \nabla \varphi^{m+1})$ still solves a linear symmetrizable 
hyperbolic system, with source terms $S^{m+1}_{\varphi}$ and $S^{m+1}_a$ already known, with 
the maximal dissipative boundary condition \iref{transcondbord}. It has then a smooth solution 
by the above mentionned theory. Then, we recover $\vp^{m+1}$ as usual by
$$ \vp^{m+1}(t,x) \equiv \int_0^t \big( S^{m+1}_\vp - 2a\, a^{m+1} - u \cdot u^{m+1}\big) (\tau,x)\ d\tau.$$
Finally, the first equation in \eqref{eqlateA+Phi+} (with $k=m+1$) is a linear ODE for $A^{m+1}$, with 
source term $F^{m+1}\in \BS_{exp}$ now determined, for which we can write down the unique explicit 
exponentially decreasing solution satisfying $\partial_{Z}A^{k}(t,y,0)= - \p_z a^k(t,y,0)$.


Consequently, we have constructed an approximate solution of \eqref{eulkort} such that
$$\left\{ \begin{array}{ll}
\displaystyle{\p_t a^\e + \nabla \vp^\e \cdot \nabla a^\e + \frac{1}{2}\, a^\e \Delta \vp^\e} & = 
\displaystyle{\e^m \big( R_{a}^{int, m}(t,x) + \e^{-1} R_{a}^{\flat,m}(t,y,z/\eps)\big)}\\ \\
\displaystyle{\p_t \vp^\e + \frac{1}{2}\, | \nabla \vp^\e|^2 + \big(a^\e \big)^2 - 1} & = 
\displaystyle{\frac{ \eps^2 }{2} \frac{\Delta a^\eps}{a^\eps}(t,x) + 
\eps^m\big( R^{int, m}_{\varphi}(t,x) + R^{\flat,m}_{\vp}(t,y,z/\eps)\big)}
\end{array}\right.
$$
where $R_{a}^{int, m}(t,x)$, $R_{\varphi}^{int,m}(t,x)$ are smooth bounded functions and
$R_{a}^{\flat,m}$, $R_{\varphi}^{\flat,m}\, \in \mathcal{S}_{exp}$. We can 
thus write the error $R^\eps$ in the GP equation as
$$ R^\eps(t,x) =  \eps^m \Big( - a^\eps( R^{int, m}_{\varphi} + R^{\flat,m}_{\varphi})
+ i ( \eps  R_{a}^{int, m} + R_{a}^{\flat,m})\Big).$$
This ends the proof of Lemma \ref{WKB+}. \b

\subsection{Validity of the WKB expansion}

\ \indent We shall now prove the stability of the WKB expansion built in Lemma \ref{WKB+}.

\begin{theo}
\label{theoGP+}
Let $\Psi^{a,m}= a^\eps e^{ i \frac{\vp^\e}{\e}}$ a WKB expansion defined on $[0, T_{m}]$ given by 
Lemma \ref{WKB+}. Then for $d \leq 3$, and $m \geq 4$ there exists a unique smooth solution 
$\Psi^\eps$ also defined on $[0, T_{m}]$ of \eqref{GP+}, \eqref{dir+}, \eqref{inf+} such that 
$ \Psi^\eps_{/t=0}= \Psi^{a,m}_{/t=0}$. Moreover, we have the estimate
$$ \e \big|\!\big|\Psi^\eps e^{-i \varphi^\eps \over \eps} - a^\eps \big|\!\big|_{H^1} 
+ \eps^3 \big|\!\big|\Psi^\eps e^{-i \frac{\vp^\e}{\e}} - a^\eps \big|\!\big|_{H^3} 
\leq C_m \eps^{m- \frac{1}{2}}, \quad \forall t \in [0, T_{m}]$$
and in particular
\beq
\label{W1infty}
\big|\!\big|\Psi^\eps e^{-i \frac{\vp^\e}{\e}} - \big(a+ \eps A^1 \big) \big|\!\big|_{W^{1, \infty}} 
\leq C_m \max \{\eps, \e^{m-\frac{7}{2}} \}.
\eeq
\end{theo}

\begin{rem}\rm 
For simplicity, we have restricted ourselves to dimension $d\leq 3$. Note however that it is possible 
to get $H^s$ estimates for every $s$. By contrast with Theorem \ref{Conver}, we emphasize that the initial 
condition in Theorem \ref{theoGP+} is exactly the WKB approximate solution $\Psi^{a,m}$. In particular, 
this initial datum has to verify some compatibility condition on the boundary.
\end{rem}

\subsubsection*{Proof.}
 As in the proof of Theorem \ref{WKBstab}, we  set 
  $$ \Psi^\eps = \Psi^{a, m} + w  \,e^{i  \varphi^\eps \over \eps  }.$$
    and we study the equation for $w$ i.e. \eqref{NLSw}. Note that we are now seeking for $w$  which
     tends to zero at infinity
     since the boundary condition at infinity is already absorbed in the WKB expansion.
    Again the first step is to get estimates for the linear  equation \eqref{NLSwlin}
     in $\Omega$ with the Neumann boundary condition
     \beq
     \label{neumw}
     \partial_{z} w(t,y,0)=0.
     \eeq
     As we can check in the proof of Lemma \ref{dtN}, in  all the integration by parts
      that are performed, the boundary terms vanish due to the Neumann
       boundary condition or the fact
        that $u^\eps_{d}(t,y,0)=0$,  and hence the proof of the $L^2$ stability  will
         be  almost the same
        as  the one in the whole space. Nevertheless, we  have to pay attantion to          the presence of boundary layer terms in the coefficients.
          At first, we note that  since $\Phi^1=0$ and $\Phi^2=0$
           in the WKB expansion, we still have that $M$
            (which is defined in Lemma \ref{dtN}) is independent of $\eps$.
            Indeed,  for  the worse term which is $\nabla (\nabla \cdot u^\eps)$, we have
           $$ \nabla(\nabla \cdot  u^\eps )  =  \partial_{ZZZ} \Phi^3 + \nabla \Delta \varphi + \mathcal{O}(\eps).$$
              
         Next,  keeping  the definitions of $R_{a}$ and $R_{\varphi}$ given
         in \eqref{Rphi}, \eqref{Ra} and by construction of the WKB expansion, we have
       \beq
       \label{Ra+} 
        |\!| R_{a}|\!|_{L^\infty} \leq C \eps^m.
        \eeq
    Nevertheless, again by construction of the WKB expansion, we only have
    $$ R_{\varphi} = R^m_{\varphi}  + {\eps^2 \over  2} { \Delta a^\eps \over a^\eps} $$
     and  due to the presence of boundary layers in $a^\eps$, we can split $R_{\varphi}$
      into
    \beq
    \label{Rphiexp}
     R_{\varphi}= \eps^2  R_{\varphi}^{int}(t,y,z) +  \eps R_{\varphi}^\flat(t, y, { z \over \eps } ) \eeq
     
    where $R_\varphi^{int}$ is smooth and bounded whereas $R_{\varphi}^\flat \in \mathcal{S}_{exp}$
     and we see that $  \eps\, |\!| R_{\varphi}^\flat |\!|_{L^\infty}= \mathcal{O}( \eps) $,  $ \eps\, |\!|\nabla R_{\varphi}^\flat|\!|_{L^\infty }
      = \mathcal{O}(1)$,  hence
       \eqref{lin0}  would be useless.
       Moreover, the fact that $R_{\varphi}^\flat$  belongs to $\mathcal{S}_{exp}$ does not
        seem to improve the estimates.
         The way to overcome this difficulty  seems to incorporate this new singular
          term into the functional.
      Let us define the operator
      $$ \mathcal{S}^\eps_{+} w  = - { \eps^2 \over 2} \Delta w + 2 (w, a^\eps) a^\eps
       + \eps R_{\varphi}^\flat w,$$
        our weighted norm in this section will be
    $$ N^\eps_{+} (w) = \int_{\Omega}  \Big( (\mathcal{S}^\eps_{+}w, w) +   K\, \eps^2 \, |w|^2 \Big) \, dx
      = { 1 \over 2  } \int_{\Omega}
       \Big( \eps^2 |\nabla w |^2 + 4 (w, a^\eps)^2 +   2 \eps R_{\varphi}^\flat \, |w|^2
        + 2 K \, \eps^2  | w|^2 \Big ) \, dx.$$
        Note that  $R_{\varphi}$ has no sign, nevertheless,  $ N^\eps_{+}(w)$  can be bounded from below
         by a weighted $H^1$ norm if $K$ is chosen sufficiently large.
          Indeed, since $R_{\varphi}^\flat$ belongs  to  $\mathcal{S}_{exp}$ we can write
       $$  2 \eps \Big|  \int_{\Omega} R_{\varphi}^\flat \, |w|^2  \, dx
        \Big| \leq C \eps \int_{\Omega}  e^{ - { \gamma z} \over \eps } |w|^2\, dx $$
         and then use 
          the one-dimensional  Sobolev inequality
        $$ | w(t, y, z) |^2 \leq C  \Big(  \int_{\mathbb{R}_{+}} |w(t,y,z ) |^2 \, dz \Big)^{ 1 \over 2 }
        \, \Big(  \int_{\mathbb{R}_{+}} | \partial_{z }w(t,y,z ) |^2 \, dz \Big)^{ 1 \over 2 }$$
        to get
    \beq
    \label{trick+}
     \eps \int_{\Omega} e^{- { \gamma z \over \eps}} |w|^2 \leq C
      \eps  |\!| w |\!|_{L^2} \, |\!| \nabla w |\!|_{L^2} \, \int_{\Omega} e^{ - { \gamma z  \over \eps } }
      \, dz
    \leq C \eps^2   \,  |\!| w |\!|_{L^2} \, |\!| \nabla w |\!|_{L^2}.
    \eeq
    In particular, we have proven that 
    \beq
    \label{trick1}
      2 \eps \Big|  \int_{\Omega} R_{\varphi}^\flat \, |w|^2  \, dx  \Big|
     \leq C \eps^2   \,  |\!| w |\!|_{L^2} \, |\!| \nabla w |\!|_{L^2}.
      \eeq
      This yields thanks to the Young inequality
     \beq
     \label{trick2}
       2 \eps \Big|  \int_{\Omega} R_{\varphi}^\flat \, |w|^2  \, dx  \Big|
     \leq  { 1 \over 2} \eps^2 |\!| \nabla w |\!|^2_{L^2}  + C \eps^2 |\!| w |\!|_{L^2}^2
     \eeq
      where $C$ is independent of $ \eps$. Consequently, if $K$ is chosen such that
       $2K>C$,
      we get
      $$ N^\eps_{+}(w) \geq  C_{0} \Big(  \eps^2 |\!|w |\!|_{H^1}^2 + \int_{\Omega} (w, a^\eps)^2 \, dx \Big), 
       \quad C_{0} >0.$$
       Note that  in this section, we have
       $$ a^\eps = a + \mathcal{O}(\eps)$$
        with $ a \geq \alpha$, this finally yields  that  $N^\eps_{+}(w)$ is equivalent
         to the weighted norm
       \beq
       \label{low+}
         N^\eps_{+}(w) \sim   \eps^2 |\!|w |\!|_{H^1}^2  + |\!|\mbox{Re } w |\!|_{L^2}^2.
         \eeq
         
        The first step in  the proof of Theorem \ref{theoGP+} is to prove the equivalent
         of Lemma \ref{dtN}. We shall prove the estimate
      \begin{eqnarray}
      \label{L2+}
\frac{d}{dt} N^\eps\big( w(t) \big)&  \leq & \ 
C  N^\eps\big( w(t) \big) 
 \\
\nonumber  & &  +   |\!| F^\eps |\!|_{L^2}^2 +  \int_{\Omega} \frac{4}{\e} (w,a^\e) (i a^\e, F^\e)
 - \int_{\Omega } (i \eps \Delta w,   F ^\eps) - \int_{\Omega}(iF^\eps, R_{\varphi}^\flat w)
\end{eqnarray} 
where $C$ independent of $\eps$.

\subsubsection*{Proof of  \eqref{L2+}}

The proof follows the same lines as the proof of  Lemma \ref{dtN}. 
 At first, since $\mathcal{S}^\eps_{+}$ is self adjoint, we have
$$
 {d \over dt} \int_{\Omega} \big( \mathcal{S}^\eps_{+} w, w \bigr) dx=
  \int_{\Omega} \Big( 2 \big( \mathcal{S}^\eps_{+} w,\partial_{t}w \big)    + 4  (w, a^\eps) (w, \partial_{t} a^\eps)
   +    2 \eps \, \partial_{t}R_{\varphi}^\flat\,  |w|^2 \Big) \,dx.
$$
Since $ \partial_{t} R_{\varphi}^\flat \in \mathcal{S}_{exp}$, we can  still use
 \eqref{trick+} to get
 $$  2 \eps \, \int_{\Omega } \partial_{t}R_{\varphi}^\flat\,  |w|^2 \leq C  N^\eps_{+}(w)$$
 Next, as in  the proof of Lemma \ref{dtN},  we use  \eqref{NLSwlin} to express $\partial_{t}w$  as
$$ \partial_{t} w = -  {i \over \eps} \mathcal{S}_{+}^\eps w  - \big( u ^\eps \cdot \nabla w +  { 1 \over 2 }w\, \nabla \cdot u^\eps \big)
 - i {  \eps^2  R_{\varphi}^{int}  \over \eps} w  - {i F^\eps \over \eps }$$
  to get
\beq
\label{estdtN1+}
2\int_{\Omega} \big(\partial_{t}w, \mathcal{S}^\eps w \big) dx
 = 2 \int_{\Omega} \Big( - \big( u^\eps \cdot \nabla w +{1 \over 2 } w \, \nabla \cdot u^\eps \big) - {i  \eps^2 R_{\varphi}^{int} \over \eps} w -
    i { F^\eps \over \eps }, \mathcal{S}^\eps_{+} w \Big) dx.
    \eeq
The only term in the right-hand side of \eqref{estdtN1+} which  is  not present in \eqref{estdtN1}, 
 is 
 $$ \mathcal{I}  =  - 2 \int_{\Omega} \Big( u^\eps \cdot \nabla w + { 1 \over 2 }w \, \nabla \cdot u^\eps, \eps R_{\varphi}^\flat w  \Big). $$
 Indeed, we have the cancellation
 $$ \quad  \int_{\Omega}( i R_{\varphi}^{int} w, R_{\varphi}^\flat w) = 0$$
  since $R_{\varphi}^{int}$ and $R_{\varphi}^\flat$ are real. To estimate $\mathcal{I}$, we note that
   we have a bound on the second term by using again \eqref{trick+}. It remains to estimate
   the first term. Integrating by parts and using that $u^\eps_{d}(t,y,0)=0$, we get
  $$    - 2 \int_{\Omega} \big( u^\eps \cdot \nabla w  , \eps R_{\varphi}^\flat w  \big)
  = \eps \int_{\Omega} \nabla \cdot\big( R_{\varphi}^\flat u^\eps \big) |w|^2
   = \int_{\Omega}  \nabla \cdot u^\eps \,  \eps R_{\varphi}^\flat |w|^2 +  \int_{\Omega}  \eps\, 
    u^\eps \cdot \nabla R_{\varphi}^\flat \, |w|^2.
   $$
   Again, the first term can be bounded thanks to \eqref{trick+}. For the second one, we first notice
    that  since $u^\eps_{d}(t,y,0) = 0$ and $R_{\varphi}^\flat \in \mathcal{S}_{exp}$,  we have
    $$  \eps \,  \big| u^\eps \cdot \nabla R_{\varphi}^\flat \big|
     \leq  C \eps \Big( |\nabla_{y}R_{\varphi}^\flat | +   | z \partial_{z} R_{ \varphi}^\flat |  \Big)
      \leq C \eps e^{ - { \gamma z \over \eps } }.$$
      This   finally yields
   $$ \mathcal{I}\leq C N^\eps_{+}(w)$$
    thanks to a new use of \eqref{trick+}.
    
    The end of the proof of \eqref{L2+} is then exactly the same as the proof of Lemma 1, since
     all the integration by parts  do not create boundary terms either because of the Neumann
      boundary condition or because $u^\eps_{d}$ vanishes on the boundary.
      
  \subsubsection*{Higher order estimates}
   
   The estimates of higher order derivatives are  more involved  than in the whole space.
   There are  two main reasons.  The first one is that    there is a new singular term $ \eps R_{\varphi}^\flat w$
    which creates bad terms when we take the derivatives of the equation.
     The second reason  is 
     that   to recover estimates on the normal derivatives, we need  to use the equation
      which gives in particular that  $ \eps^2 \partial_{z}^2$ behaves like
       $ \eps \partial_{t}$ and $\eps \nabla$. This anisotropy in the weights 
        does not seem to allow to construct  high order functionals like $N^\eps_{s}(w)$
         which allows  to get $H^s$ estimates without  additional loss of  $ \eps $.
       Let us use the notation 
       $$\Lambda= (\Lambda_{0}, \cdots \Lambda_{d}) = (\partial_{t}, \nabla_{y}, p(z) \partial_{z})^t$$
        where the weight $p(z)$
         is given by $p(z)= z/(1+z).$ Note that we can apply $\Lambda$
          to the equation since $\Lambda w$ still satisfies the Neumann
           boundary condition.  The use of $\Lambda$ is
           classical in hyperbolic  characteristic initial boundary value problems
            (see \cite{Rauch} for example)
         The weighted norm that we shall estimate is
       $$ Y^\eps_{+}(w) = N^\eps_{+} (w) + N^\eps_{+}( \eps \Lambda w).$$
       In dimension $d \leq 3$, this is sufficient to get the nonlinear stability.
        We shall see in the proof why the use of $\Lambda_{d}$ is necessary
        
       We shall prove that
     \beq
     \label{H2+}
      {d \over dt } Y^\eps_{+}(w)
       \leq  C\Big(  Y^\eps_{+}(w)    
           + X^\eps(F^\eps) + X^\eps( \eps \Lambda F^\eps) 
          \Big)
     \eeq
      for some $C>0$ independent of $\eps$ where
      we have set
      $$ X^\eps(F)=  |\!|F|\!|_{H^1}^2 + { |\!| F |\!|_{L^2}^2 \over \eps }
       + {  |\!| \mbox{Im } F |\!|_{L^2}^2 \over \eps^2}.$$
      
    \subsubsection*{Proof of \eqref{H2+}}
    As a preliminary, we shall rewrite \eqref{L2+} in a more convenient form.  We can
     use that $a^\eps = a+ \mathcal{O}(\eps)$ with $a $ real, perform an integration
      by parts and use \eqref{trick+}  to get from \eqref{L2+} that
   \beq
   \label{N2+2}
    \frac{d}{dt} N^\eps\big( w(t) \big)&  \leq & \ 
C  N^\eps\big( w(t) \big) +  X^\eps (F^\eps )
\eeq
where
$$ X^\eps (F^\eps) = |\!|F^\eps |\!|_{H^1}^2 + { |\!|F |\!|_{L^2}^2 \over \eps  }
 + { |\!| \mbox{Im } F^\eps |\!|_{L^2}^2 \over \eps^2}.
 $$
 
 To prove \eqref{H2+}, we  start with the estimate of $N^\eps_{+}(\eps \partial_{t} w )$.
  When we apply $\eps \partial_{t}$ to \eqref{NLSwlin}, we find
 \beq
 \label{NLSwlindt} \Big(  i \eps \partial_{t} + \mathcal{L}^\eps \Big) \eps \partial_{t}w = R_{\varphi} \, \eps \partial_{t} w 
  + \eps \partial_{t} F^\eps + \mathcal{C}
  \eeq
  where the commutator $\mathcal{C}$ can be splitted into
  \beq
  \label{Cdef+} \mathcal{C}= \mathcal{C}_{1} + \mathcal{C}_{2}+ \mathcal{C}_{3}
  \eeq
   with 
   \begin{eqnarray*}
   & & \mathcal{C}_{1 }= \eps  \partial_{t}R_{\varphi} w, \\
   & & \mathcal{C}_{2}  = 2  \eps \Big( ( \partial_{t}a^\eps, w) a^\eps + ( a^\eps, w) \partial_{t} a^\eps \Big), \\
   & & \mathcal{C}_{3}=  - i \eps^2 \Big( \partial_{t} u^\eps \cdot \nabla w +  { 1 \over 2 }\partial_{t}(\nabla \cdot u^\eps)\, w \Big).
   \end{eqnarray*}
  Consequently, we can apply \eqref{N2+2} to \eqref{NLSwlindt} with the new source
   term $\eps \partial_{t}F^\eps + \mathcal{C}$ to get
   \beq
   \label{dt+1}
    \frac{d}{dt} N^\eps_{+}\big( \eps \partial_{t} w(t) \big)&  \leq & \ 
C  N^\eps_{+}\big(  \eps \partial_{t} w(t) \big) +  X^\eps ( \eps \partial_{t} F^\eps ) + X^\eps (\mathcal{C}).
\eeq
Thus it remains to estimate $X^\eps (\mathcal{C})$.
  Let us begin with $X^\eps( \mathcal{C}_{1}).$ Thanks to the expansion \eqref{Rphiexp}, we easily get
 \begin{eqnarray}
 \nonumber  X^\eps(\mathcal{C}_{1})&  \lesssim&  N^\eps_{+}(w) +  \int_{\Omega}  \Big( |\partial_{t}R_{\varphi}^\flat |^2\,  
   \eps^4  |w|^2 +
 \eps^4  |\nabla w |^2 )
  + \eps^4 |\nabla \partial_{t}R_{\varphi}^\flat |^2 \,  |w|^2  \Big)\,dx. \\
 \label{C1+} & \lesssim &    N^\eps_{+}(w) . \end{eqnarray}
  Note that we could have a better estimate  by using that $R_{ \varphi}^\flat \in \mathcal{S}_{exp}$
   and  \eqref{trick+}. Next, we turn to the estimate of
    $X^\eps(\mathcal{C}_{2})$. By using that  $ a^\eps = a+ \mathcal{O}(\eps)$
     with $a$ real, we find
    \beq
    \label{C2+}
    X^\eps(\mathcal{C}_2)  \lesssim  N^\eps_{+}(w) + \eps   |\!|\mbox{Re }w |\!|_{L^2}^2 
    + \eps^2 |\!|\nabla w |\!|_{L^2}^2 \lesssim N^\eps(w).
    \eeq
    Note that the above estimate was sharp. This  is for the estimate of this commutator $\mathcal{C}_{2}$
     that  we had to  chose the weight  $  \eps$ in front of the time derivative. Finally, we estimate
      $X^\eps(\mathcal{C}_{3})$. For the estimate of $\mathcal{C}_{3}$, we  use that
       $\partial_{t}u^\eps_{d}$ vanishes on the boundary which implies that
       $$ | \partial_{t}u^\eps_{d} | \lesssim p(z).$$
        Thanks to this remark, we find
   \begin{eqnarray}
 & &   \label{C3+} X^\eps (\mathcal{C}_3) \lesssim  N^\eps_{+}(w) +  \eps^4 |\!| \nabla  \Lambda w |\!|_{L^2}^2
   \lesssim Y^\eps_{+}(w).
   \end{eqnarray}
   Note that this is  for the control of this commutator that we are obliged to add the vector field
    $p(z) \partial_{z}$ in the definition of  the functional space.
     Consequently, the combination of \eqref{dt+1}, \eqref{C1+}, \eqref{C2+} and \eqref{C3+} gives
   \beq
   \label{dt+2}
    \frac{d}{dt} N^\eps_{+}\big( \eps \partial_{t} w(t) \big)&  \lesssim  &  
  Y^\eps_{+}\big(   w(t) \big) +  X^\eps ( \eps \partial_{t} F^\eps ).
  \eeq
  The estimate of $\eps \nabla_{y}w$ follows exactly the same lines, and we also find
   \beq
   \label{dt+3}
    \frac{d}{dt} N^\eps_{+}\big( \eps \nabla_{y} w(t) \big)&  \lesssim  &  
  Y^\eps_{+}\big(  w(t) \big) +  X^\eps ( \eps \nabla_{y} F^\eps ).
  \eeq
  The estimate of $ \eps \Lambda_{d} w= \eps p(z)  \partial_{z}w$ requires  some additional
   work since the vector field $\Lambda_{d}$ does not commute with the Laplacian.
    By applying $ \eps \Lambda_{d}$ to \eqref{NLSwlin}, we get
    \beq
 \label{NLSwlindz} \Big(  i \eps \partial_{t} + \mathcal{L}^\eps \Big) \eps \Lambda_{d}w = R_{\varphi} \eps \Lambda_{d} w 
  + \eps \Lambda_{d} F^\eps + \mathcal{C} + \mathcal{C}_{4}
  \eeq
where $\mathcal{C}$ is defined as in \eqref{Cdef+} above with $\partial_{t}$ replaced by $\Lambda_{d}$ and
 $\mathcal{C}_{4}$ is given by
 $$ \mathcal{C}_{4} =  - { \eps^3 \over 2}[ \Lambda_{d}, \Delta ] =  - { \eps^3 \over 2}
  \Big( 2 \partial_{z}p \, \partial_{zz}w +  \partial_{zz}p \, \partial_{z}w \Big).
  $$ 
 Next, we can apply  \eqref{L2+}  to get   
\begin{eqnarray*} 
 \frac{d}{dt} N^\eps_{+}\big( \eps \Lambda_{d} w(t) \big)  & \lesssim &     
  N^\eps_{+}\big( \eps \Lambda_{d}  w(t) \big) +  X^\eps ( \eps \Lambda_{d} F^\eps )
   + X^\eps(\mathcal{C}) + |\!| \mathcal{C}_{4}|\!|_{H^1}^2 \\
  & &  + { 4 \over \eps} \int_{\Omega} ( \eps \Lambda_{d}w, a^\eps) (i a^\eps,  
  \mathcal{C}_{4})  - \int_{\Omega}\big(i  \mathcal{C}_{4}, R^\flat_{\varphi}\, \eps
    \Lambda_{d} w \big)
\end{eqnarray*}
Since  one can easily check that $X^\eps(\mathcal{C})$ still satisfies the bounds
   \eqref{C1+}, \eqref{C2+}, \eqref{C3+}, we obtain
 \begin{eqnarray*} 
 \frac{d}{dt} N^\eps_{+}\big( \eps \Lambda_{d} w(t) \big)  & \lesssim &  
  Y^\eps_{+}(w)   +  X^\eps ( \eps \Lambda_{d} F^\eps ) 
 + |\!| \mathcal{C}_{4}|\!|_{H^1}^2 \\
  & &  + { 4 \over \eps} \int_{\Omega} ( \eps \Lambda_{d}w, a^\eps) (i a^\eps,  
  \mathcal{C}_{4})  - \int_{\Omega}\big(i \mathcal{C}_{4}, R^\flat_{\varphi}\, \eps
    \Lambda_{d} w \big).
\end{eqnarray*}
Next, we note that 
 $$  |\!| \mathcal{C}_4 |\!|_{H^1}^2  \lesssim \eps^6 |\!|w|\!|_{H^3}^2$$
  and that
 \begin{eqnarray*}
 { 4 \over \eps} \Big|  \int_{\Omega} ( \eps \Lambda_{d}w, a^\eps) (i a^\eps,  
  \mathcal{C}_{4}) \Big| \lesssim  {4 \over \eps} \int_{\Omega}
     \eps | \partial_{z} w | \, | p(z) \mathcal{C}_{4}  |
     & \lesssim  &  \eps^2 N^\eps_{+}(w)^{ 1 \over 2} \Big( |\!|p \partial_{zz} w |\!|_{L^2}
       + |\!| \partial_{z}w |\!|_{L^2} \Big) \\
       & \lesssim& N^\eps_{+}(w)^{ 1 \over 2} Y^\eps_{+}(w)^{ 1 \over 2}.
 \end{eqnarray*}
  In a similar way, we also get 
$$ \Big|Ê \int_{\Omega}\big(i \mathcal{C}_{4}, R^\flat_{\varphi}\, \eps
    \Lambda_{d} w  \Big| \lesssim  \eps |\!| \partial_{z} w  |\!  |_{L^2} \,  |\!| p\,  \mathcal{C}_{4}|\!|_{L^2}
     \lesssim  Y^\eps_{+}(w).$$
 Consequently, we have proven that
 \beq
 \label{lambdad+}  \frac{d}{dt} N^\eps_{+}\big( \eps \Lambda_{d} w(t) \big)   \lesssim   
  Y^\eps_{+}(w)   +  X^\eps ( \eps \Lambda_{d} F^\eps )  + \eps^6 |\!| w |\!|_{H^3}^2.\eeq
   To conclude,  it remains to estimate $ \eps^6 |\!| w |\!|_{H^3}^2$.
    As usual, this is done thanks to the equation \eqref{NLSwlin} and the standard regularity result
     for elliptic equations. We  rewrite \eqref{NLSwlin} as the equation
     \beq
     \label{elliptic+}
      \eps^2 \Delta w = G^\eps, \quad \partial_{z}w(t,y,0)=0
      \eeq
       where the source term enjoys the estimates
       \begin{eqnarray*}
        & & |\!| G^\eps |\!|_{L^2}^2 \lesssim  \eps^2 |\!| \Lambda w|\!|_{L^2}^2  + |\!|w|\!|_{L^2}^2
         + |\!| F^\eps |\!|_{L^2}^2 , \\
        & & |\!| \nabla G^\eps |\!|_{L^2}^2 \lesssim \eps^2  |\!| \nabla \Lambda w |\!|_{L^2}^2
         + |\!| w |\!|_{H^1}^2 + |\!| \nabla F^\eps |\!|_{L^2}^2.
         \end{eqnarray*}
         Consequently, we  get from \eqref{elliptic+} by standard elliptic regularity that
         \beq
         \label{H3+} \eps^6 |\!| w |\!|_{H^3}^2 \lesssim Y^\eps_{+}(w) +  |\!| F^\eps |\!|_{H^1}^2.
         \eeq
          By replacing this last estimate in \eqref{lambdad+},  we finally  obtain
       \beq
       \label{zdz} \frac{d}{dt} N^\eps_{+}\big( \eps \Lambda_{d} w(t) \big)   \lesssim   
  Y^\eps_{+}(w)   +  X^\eps ( \eps \Lambda_{d} F^\eps ) +  |\!|F^\eps |\!|_{H^1}^2.
  \eeq
  
   To conclude, it suffices to sum the estimates \eqref{N2+2}, \eqref{dt+2}, \eqref{dt+3}
    and \eqref{zdz} to get \eqref{H2+}. 
\bigskip

The estimate \eqref{H2+} is sufficient to prove the nonlinear stability stated
 in Theorem \ref{theoGP+} for $d\leq 3$. Nevertheless, it is
  possible  to prove by induction that
   for every $s$,
   $$ { d \over dt } \Big( \sum_{m \leq s} N^\eps_{+}( (\eps\Lambda)^m w )  \Big)
    \lesssim   \sum_{m \leq s} \Big( X^\eps( (\eps \Lambda)^m F^\eps) + N^\eps_{+}( (\eps\Lambda)^m w )  \Big) .
    $$

    \subsection*{Nonlinear stability}
  Thanks to \eqref{H2+} and Gronwall inequality, we get for $0 \leq T\leq T_{m}$
  $$ \sup_{[0, T]} Y^\eps_{+}(w) \lesssim Y^\eps_{+}(0) +  Te^{\gamma T} \sup_{[0,T]}\Big(
   X^\eps( F^\eps) + X^\eps( \eps \Lambda F^\eps) \Big)$$
   for some $\gamma>0$ independent of $\eps$.
   Consequently, we can combine this last estimate with \eqref{H3+} to get
   \beq
   \label{apriori+}
   \sup_{[0, T]}  Z^\eps_{+}(w) \lesssim  Y^\eps_{+}(0) +   (1+T)e^{\gamma T} \sup_{[0,T]}\Big(
   X^\eps( F^\eps) + X^\eps( \eps \Lambda F^\eps) \Big)
   \eeq
   with
   $$ Z^\eps_{+}(w) =  Y^\eps_{+}(w) + \eps^6 |\!| w |\!|_{H^3}^2.$$
   Thanks to this a priori estimate, one can easily prove by standard fixed point argument
    the   existence  of a unique
    solution of \eqref{NLSw} on some interval of time $[0, T^\eps]\subset[0, T_{m}]$
      such that $Z^\eps_{+}(w)$ remains finite.
      
      By using  that $w_{/t=0}=0$ and the equation to compute the time derivative, we find
      $$ Y^\eps_{+}(w(0)) \leq C_{0} \eps^{2m}.$$
       Moreover, using that $F^\eps= \eps^m R^\eps + Q^\eps,$
        we have thanks to \eqref{formR+} that
       $$  (1+T_{m})e^{\gamma T_{m}}\sup_{[0, T_{m}]}\Big(
        X^\eps_{+}(R^\eps) + X^\eps_{+}(\Lambda  R^\eps) \Big) \leq C(T_{m}) \eps^{2m-1}$$
        where $T_{m}$ is the existence time of the approximate solution
         given by Lemma \ref{WKB+}. We  can thus  fix $R>C_{0}+ C(T_{m})$ and define
          $\tau^\eps$ the maximal time for which the solution of \eqref{NLSw} satisfies
          $ Z^\eps_{+}(w(t)) \leq R \eps^{2m - 1}.$ As in the proof of Theorem \ref{WKBstab}, we shall
           prove that for $\eps$ sufficiently small, we have $\tau^\eps =  T_{m}$.
            Thanks to \eqref{apriori+},  we have for every $T$, $T<\tau^\eps$, 
        \beq
        \label{apriori2+}
        \sup_{[0, T]}  Z^\eps_{+}(w) \lesssim  R \eps^{ 2m - 1} +   (1+T)e^{\gamma T} \sup_{[0,T]}\Big(
   X^\eps( Q^\eps) + X^\eps( \eps \Lambda Q^\eps) \Big).
   \eeq
   Here, the expression of $Q^\eps(w)$ is given by
   $$ Q^\eps(w) =  a^\eps |w |^2 + 2 (a^\eps, w) w +w |w|^2.$$
    To conclude, we need to bound the right hand side of \eqref{apriori2+}. To estimate the nonlinear
     term, we  use that for $d \leq 3$, we have
     \beq
     \label{sob+}
      |\!| w|\!|_{L^\infty}^2 \lesssim |\!| \nabla^2 w |\!|\, |\!| w | \! |_{H^1},       \eeq
      which gives
     \beq
     \label{winf+}
       |\!| w|\!|_{L^\infty}^2 \lesssim { Z^\eps_{+}(w)  \over \eps^4}\lesssim \eps^{2m-5}\,
       \quad \forall t \in [0, \tau^\eps). 
      \eeq
    We shall take $m$ such that $2m>5$ in order to get 
   $  |\!| w|\!|_{L^\infty} \leq 1$  for $t \in [0, \tau^\eps).$  
      This implies
     $$ |\!| Q^\eps |\!|_{H^1}^2 \lesssim \big(  |\!| w |\!|_{L^\infty}^2
      + |\!| w |\!|_{L^\infty}^4 \big)  \, |\!| w |\!|_{H^1}^2 \lesssim  { Z^\eps_{+}(w)^2
      \over \eps^6}.$$
      Next, since $H^1(\mathbb{R}^d) \subset L^4$ for $ d \leq 3$, we also have
     $$ { |\!| Q^\eps |\!|_{L^2}^2  \over \eps^2  } \lesssim { |\!| w |\!|_{H^1}^4 \over  \eps ^2 } 
      ( 1  +  |\!| w |\!|_{L^\infty}^2 )   \lesssim 
      { Z^\eps_{+}(w)^2 \over  \eps^6}  .$$
      Consequently, we have already proven that
      \beq
      \label{Q1+}
      X^\eps (Q^\eps)  \lesssim  { Z^\eps_{+}(w)^2 \over  \eps^6}.
      \eeq
      Next, we evaluate $X^\eps( \eps \Lambda Q^\eps)$.
       At first, we write
       \begin{eqnarray*}
        \eps^2  |\!| \Lambda Q^\eps |\!|_{H^1}^2 \lesssim 
             \eps^2  |\!| \Lambda w |\!|_{H^1}^2 \big(  |\!| w |\!|_{L^\infty}^2 + 
              |\!| w | \!|_{L^\infty}^4 \big)
               +  \eps^2 |\!| \Lambda w |\!|_{L^4}^2 \, |\!| \nabla w |\!|_{L^4}^2 \big( 1 + |\!| w |\!|_{L^\infty}^2
                ) 
        \end{eqnarray*}
         and by using  for $d \leq 3$, the Sobolev embedding
          $H^1 \subset L^4 $ and the Gagliardo-Nirenberg inequality
         $$  |\!| \nabla f | \!|_{L^4}^2 \lesssim |\!|f |\!|_{H^1}^{ 1 \over 2 }\, |\!| \nabla^2
f |\!|_{L^2}^{ 3 \over 2 } , $$
we get  for $0 \leq t \leq \tau^\eps$: 
$$   \eps^2  |\!| \Lambda Q^\eps |\!|_{H^1}^2
 \lesssim  { Z^\eps_{+}(w)^2 \over \eps^4 } +  \eps^2   |\!| \nabla w |\!|_{H^1}^2
 |\!|w |\!|_{H^1}^{ 1 \over 2 }\, |\!| \nabla^2
w |\!|_{L^2}^{ 3 \over 2 }  \lesssim    { Z^\eps_{+}(w)^2 \over \eps^6}.
$$
 Finally, by similar arguments, we also have
$$ {  |\!| \eps \Lambda Q^\eps |\!|_{L^2}^2 \over \eps^2 } \lesssim |\!|  \Lambda w |\!|
_{L^4}^{ 1 \over 2 }\, |\!|w |\!|_{L^4}^{ 1 \over 2}
 \lesssim  |\!|  \Lambda w |\!|
_{H^1}^{  2 }\, |\!|w |\!|_{H^1}^{ 2}
 \lesssim { Z^\eps_{+}(w)^2  \over  \eps^6}.$$
 We have thus proven that
 \beq
 \label{Q2+}
  X^\eps( \eps \Lambda Q^\eps) \lesssim { Z^\eps_{+}(w)^2  \over  \eps^6}.
  \eeq 
  Consequently, by using \eqref{apriori2+}, \eqref{Q1+}, \eqref{Q2+}, 
  we get 
 $$  \sup_{[0, T]}  Z^\eps_{+}(w) \lesssim  R \eps^{ 2m - 1} +  C(R)  (1+T)e^{\gamma T} \sup_{[0,T]}  { Z^\eps_{+}(w)^2  \over  \eps^6}.$$
  By choosing $m \geq 4$, this allows to  get that $\tau^\eps =  T_{m}$
   for $\eps$ sufficiently small
   and that
   $$  \sup_{[0, T_{m}]}  Z^\eps_{+}(w) \leq  C  \eps^{ 2m - 1}.$$
   Finally, the estimate  \eqref{W1infty} follows by Sobolev embedding.
    This ends the proof of Theorem \ref{theoGP+}.

 \subsection*{Aknowledgement} We thank Remi Carles for usefull
  comments about this work.

\appendix

\section{A Lemma about composition in Sobolev spaces}

\ \indent During the proof of Lemma \ref{estirhs}, we have used a result about 
composition in Sobolev spaces. This result is very standard when $h$ does not depend 
on $x$ (see, for instance, \cite{Taylor}).

\begin{lem}
\label{composition}
Let $R>0$, $s \in \N$ and $h = h (x,w)\in \BC^{s+1} (\R^d \times \R^2, \R)$, 
satisfying $h(x,0)=0$ for all $x\in \R^d$. Assume moreover
$$ A \equiv \sup \big\{ \big|\! \big| \p_x^\alpha \p_w^\beta h \big|\! \big|_{L^\ii(\R^d \times B_R)}, \ 
\alpha \in \N^d,\ \beta \in \N^2, \ |\alpha| + |\beta| \leq s +1 \big\} < + \ii. $$
Then, there exists $C$, depending only on $A$, $s$ and $R$, such that, for any 
$w\in H^s(\R^d)$ satisfying $|w|_{L^\ii(\R^d)} \leq R$, we have $h\big(x,w(x) \big) \in H^s(\R^d)$ 
and
$$ \big|\! \big| h\big(x,w(x) \big) \big|\! \big|_{H^s} \leq C \big|\! \big| w \big|\! \big|_{H^s}.$$
\end{lem}

\noindent {\it Proof.} The proof is by induction on $s\in \N$ and relies on 
the Gagliardo-Nirenberg inequality. If $s=0$, it suffices to notice that since 
$h(x,0)=0$, then for $w \in B_R$,
$$ | h(x,w) | \leq A |w|.$$
Assume then the result for $s-1\in \N$. Let $\mu \in \N^d$ with $|\mu | =s$. 
One has easily
$$ \p^\mu \big( h(x,w(x) ) \big) = 
\sum * \big( \p_x^\alpha \p_w^{\beta + \gamma} h \big) \big( x, w(x) \big) 
\big( \p^\beta w_1\big)^p \big( \p^\gamma w_2 \big)^q,$$
where $\alpha\in \N^d$, $\alpha \leq \mu$, $\beta$, $\gamma \in \N^2$, 
$p$, $q \in \N^*$ depend on $\beta$ and $\gamma$, $|\alpha| + p |\beta| + q |\gamma| = s$, 
and $*$ is a coefficient depending only on $\mu$, $\alpha$, $\beta$ and $\gamma$. 
Furthermore, since $w\in H^s \cap L^\ii$, the Gagliardo-Nirenberg inequality yields, 
for $1 \leq k \leq s$,
$$ \big|\! \big| w \big|\! \big|_{W^{k,\frac{2s}{k}}} \leq C_{k,s} 
\big|\! \big| w \big|\! \big|_{H^s}^{\frac{k}{s}} 
\big|\! \big| w \big|\! \big|_{L^\ii}^{1-\frac{k}{s}}.$$
As a consequence, by interpolation, if $w\in H^s \cap L^\ii$ and 
$\big|\! \big| w \big|\! \big|_{L^\ii} \leq R$, then for $\gamma \in \N^d$, 
$|\gamma |\leq s$, and $2 \leq p \leq \frac{2s}{|\gamma|}$,
$$ \big|\! \big| \p^\gamma w \big|\! \big|_{L^p} \leq C_{s,p,R} 
\big|\! \big| w \big|\! \big|_{H^s}^{\frac{2}{p}}.$$
Therefore, in view of $|\alpha| + p |\beta| + q |\gamma| = s$, 
by H\"older inequality, we can estimate the terms in $\p^\mu \big( h(x,w(x) ) \big)$ 
for which $\alpha \not = \mu$ (thus $|\alpha| < s$) as
$$\big|\! \big| 
\big( \p_x^\alpha \p_w^{\beta +\gamma} h \big) \big(x,w(x) \big) 
\big( \p^\beta w_1 \big)^p \big(\p^\gamma w_2\big)^q \big|\! \big|_{L^2} 
\leq A \big|\! \big| \p^\beta w_1 \big|\! \big|_{L^{2\frac{s-|\alpha|}{|\beta|}}}^p 
\big|\! \big| \p^\gamma w_2 \big|\! \big|_{L^{2\frac{s-|\alpha|}{|\gamma|}}}^q 
\leq C_{s,p,R} A\big|\! \big| w \big|\! \big|_{H^s}.$$
For the term for which $\alpha = \mu$, we note that since $h(x,0)=0$ for $x\in \R^d$, 
then $ (\p^\alpha_x h) (x,0)=0$ for any $x\in \R^d$, so that if $w\in B_R \subset \R^2$,
$$ \big| (\p_x^\alpha h)(x,w) \big| \leq A |w|,$$
which implies
$$ \big|\! \big| (\p^\alpha_{x} h) \big( x,w(x) \big) \big|\! \big|_{L^2} 
\leq A \big|\! \big| w \big|\! \big|_{L^2} \leq A \big|\! \big| w \big|\! \big|_{H^s}.$$
Combining these two estimates gives
$$ \big|\! \big| \p^\mu \big( h(x,w(x)) \big) \big|\! \big|_{L^2} 
\leq C_{s,p,R} A \big|\! \big| w \big|\! \big|_{H^s}$$
and the proof of the Lemma is complete. \b

\newpage


\end{document}